\documentclass{amsart}
\usepackage{amsmath}
\usepackage{amssymb}
\usepackage{amscd}

\input tcilatex
\theoremstyle{definition}
\theoremstyle{remark}
\numberwithin{equation}{section}

\begin{document}
\title{Quantum Stochastic Semigroups and Their Generators}
\author{V.\thinspace P.\thinspace Belavkin}
\address{Mathematics Department, University of Nottingham, NG7 2RD, UK.}
\email{vpb@maths.nott.ac.uk}
\date{November 10, 1996}
\thanks{Published in: \textit{Irreversibility and Causality} 82 --109.
Springer Verlag, Lecture Notes in Physics, Berlin, 1998.}
\maketitle

\begin{abstract}
A rigged Hilbert space characterisation of the unbounded generators of
quantum completely positive (CP) stochastic semigroups is given. The general
form and the dilation of the stochastic completely dissipative (CD) equation
over the algebra $\mathcal{L}\left( \mathcal{H}\right) $ is described, as
well as the unitary quantum stochastic dilation of the subfiltering and
contractive flows with unbounded generators is constructed.
\end{abstract}

\section*{Introduction}

Quantum stochastic dynamics gives beautiful solvable models for the
interaction of a quantum system with the quantum noise, which is produced by
a heat bath, measurement apparatus, or any other environment with infinite
number of freedom. It can be defined by a weakly continuous evolution
semigroup on a rigged Hilbert space with a special unbounded form-generator,
corresponding to the singular boundary-type interaction. The Heisenberg
picture of such interaction is described by the quantum stochastic Langevin
equation.

In quantum theory of open systems there is a well known Lindblad's form \cite%
{19} of quantum Markovian master equation, satisfied by the one-parameter
semigroup of completely positive (CP) maps over the algebra $\mathcal{L}%
\left( \mathcal{H}\right) $ of bounded operators on the system Hilbert space 
$\mathcal{H}$. This is nonstochastical equation, which can be obtained by
averaging stochastic Langevin equation for quantum flow \cite{EvH} over the
driving noises, represented in a Fock space $\mathcal{F}$. On the other hand
the quantum EH-flow corresponds to the interaction representation for a one
parametric group of dynamical authomorphisms over $\mathcal{L}\left( 
\mathcal{H}\otimes \mathcal{F}\right) $, which are obviously completely
positive due to *-multiplicativity of these representations. The
authomorphisms (representations) give the examples of pure, i.e. extreme
point CP maps, but among the extreme points of the convex cone of all CP
maps over $\mathcal{L}\left( \mathcal{H}\right) $ there are not only the
representations. This means a possibility to construct the stochastic
representations of dynamical CP semigroups as averagings of pure, i.e.
non-mixing irreversible quantum stochastic CP dynamics, which can not be
driven by a Langevin equation. Such irreversible dynamics, corresponding to
the interaction representation for a dynamical CP semigroup over $\mathcal{L}%
\left( \mathcal{H}\otimes \mathcal{F}\right) $, are described by quantum
stochastic flows of CP maps, which should satisfy a generalized form of
Lindblad equation with quantum stochastic unbounded generators.

The examples of such dynamics having recently been found many physical
applications, will be considered in the first section. The rest of the paper
will be devoted to the mathematical derivation of the general structure for
the unbounded generators of the dynamical CP semigroups, corresponding to
the quantum stochastic CP flows over $\mathcal{L}\left( \mathcal{H}\right) $
with the noises represented in $\mathcal{F}$. The results of the paper not
only generalize the Evans-Hudson (EH) flows \cite{EvH} from the
representations to the general CP maps, but also prove the existence of the
homomorphic dilations for the subfiltering and contractive CP flows. This
gives the subfiltering CP flows as conditional expectations of EH flows,
generalizing the similar representation for contractive CP semigroups. Here
in the introduction we would like to outline the generalized structure of
the generators on the formal level.

As was proved in \cite{Bcs}, every stationary quantum stochastic processes $%
t\in \mathbb{R}_{+}\mapsto \Lambda \left( t,a\right) $ parametrized by $a\in 
\mathfrak{a}$ with $\Lambda \left( 0,a\right) =0$ and independent increments 
$\mathrm{d}\Lambda \left( t,a\right) =\Lambda \left( t+\mathrm{d}t,a\right)
-\Lambda \left( t,a\right) $, forming an It\^{o} $\star $-algebra 
\begin{equation}
\mathrm{d}\Lambda \left( a\right) ^{*}\mathrm{d}\Lambda \left( a\right) =%
\mathrm{d}\Lambda \left( a^{\star }a\right) ,\quad \sum \lambda _i\mathrm{d}%
\Lambda \left( a_i\right) =\mathrm{d}\Lambda \left( \sum \lambda
_ia_i\right) ,\quad \mathrm{d}\Lambda \left( a\right) ^{*}=\mathrm{d}\Lambda
\left( a^{\star }\right) ,  \label{0.1}
\end{equation}
can be represented in the Fock space $\mathfrak{F}$ over the space of $%
\mathcal{E}$ -valued square-integrable functions on $\mathbb{R}_{+}$ as $%
\Lambda \left( t,a\right) =a_\nu ^\mu \Lambda _\mu ^\nu \left( t\right) $
with respect to the vacuum state $\delta _\emptyset \in \mathfrak{F}$. Here 
\begin{equation}
a_\nu ^\mu \Lambda _\mu ^\nu \left( t\right) =a_{\bullet }^{\bullet }\Lambda
_{\bullet }^{\bullet }\left( t\right) +a_{+}^{\bullet }\Lambda _{\bullet
}^{+}\left( t\right) +a_{\bullet }^{-}\Lambda _{-}^{\bullet }\left( t\right)
+a_{+}^{-}\Lambda _{-}^{+}\left( t\right) ,  \label{0.3}
\end{equation}
is the canonical decomposition of $\Lambda $ into the exchange $\Lambda
_{\bullet }^{\bullet }$, creation $\Lambda _{\bullet }^{+}$, annihilation $%
\Lambda _{-}^{\bullet }$ and preservation (time) $\Lambda _{-}^{+}=t\mathrm{I%
}$ processes of quantum stochastic calculus \cite{Par}, \cite{Mey} having
the mean values $\left\langle \Lambda _\mu ^\nu \left( t\right)
\right\rangle =t\delta _{+}^\nu \delta _\mu ^{-}$ with respect to the vacuum
state in $\mathfrak{F}$, and $\mathcal{E}$ is a pre-Hilbert space of the
quantum noise in $\mathfrak{F}$. Thus the parametrizing algebra $\mathfrak{a}
$ can be always identified with a $\star $-subalgebra of the algebra $%
\mathcal{Q}\left( \mathcal{E}\right) $ of all quadruples $\boldsymbol{a}%
=\left( a_\nu ^\mu \right) _{\nu =+,\bullet }^{\mu =-,\bullet }$, where $%
a_\nu ^\mu :\mathcal{E}_\nu \rightarrow \mathcal{E}_\mu $ are the linear
operators on $\mathcal{E}_{\bullet }=\mathcal{E},\mathcal{E}_{+}=\mathbb{C=}%
\mathcal{E}_{-}$, having the adjoints $a_\nu ^{\mu *}\mathcal{E}_\mu
\subseteq \mathcal{E}_\nu $, with the Hudson--Parthasarathy (HP)
multiplication table \cite{16} 
\begin{equation}
\boldsymbol{a}\bullet \boldsymbol{b}=\left( a_{\bullet }^\mu b_\nu ^{\bullet
}\right) _{\nu =+,\bullet }^{\mu =-,\bullet },  \label{0.4}
\end{equation}
and the involution $a_{-\nu }^{\star \mu }=a_{-\mu }^{\nu *}$, where $-(-)=+$%
, $-\bullet =\bullet $, $-(+)=-$.

The stochastic differential of a CP flow $\phi =\left( \phi _t\right) _{t>0}$
over an operator algebra $\mathcal{B}\subseteq \mathcal{L}\left( \mathcal{H}%
\right) $ is written in terms of the quantum canonical differentials as $%
\mathrm{d}\phi =\phi \circ \lambda _\nu ^\mu \mathrm{d}\Lambda _\mu ^\nu $
with $\phi _0=\imath $ at $t=0$, where $\imath \left( B\right) =B$ is the
identical representation of $\mathcal{B}$. The main result of this paper is
the derivation and the dilation of the linear quantum stochastic evolution
equation 
\begin{equation*}
\mathrm{d}\phi _t\left( B\right) +\phi _t\left( K^{*}B+BK-L^{*}\jmath \left(
B\right) L\right) \mathrm{d}t=\phi _t\left( L^{\bullet }\jmath \left(
B\right) L_{\bullet }-B\otimes \delta _{\bullet }^{\bullet }\right) \mathrm{d%
}\Lambda _{\bullet }^{\bullet }
\end{equation*}
\begin{equation}
+\phi _t\left( L^{\bullet }\jmath \left( B\right) L-K^{\bullet }B\right) 
\mathrm{d}\Lambda _{\bullet }^{+}+\phi _t\left( L^{*}\jmath \left( B\right)
L_{\bullet }-BK_{\bullet }\right) \mathrm{d}\Lambda _{-}^{\bullet },
\label{0.5}
\end{equation}
where $\jmath $ is an operator representation of $\mathcal{B}$, $\delta
_{\bullet }^{\bullet }$ is the identity operator in $\mathcal{E}$, and the
operator $K$ satisfies the dissipativity condition $K+K^{\dagger }\geq
L^{*}L $ with the Hamiltonian part $H=\func{Im}K$. This differential form
for the CP flows was discovered in \cite{15} as the general completely
dissipative (CD) structure of the bounded quantum stochastic generators $%
\lambda _\nu ^\mu :\mathcal{B}\rightarrow \mathcal{B}$ over a von Neumann
algebra $\mathcal{B}$ even in the nonlinear case. In the matrix form $%
\boldsymbol{\lambda }=\left( \lambda _\nu ^\mu \right) _{\nu =+,\bullet
}^{\mu =-,\bullet }$ this can be written similar to the Lindblad form for
the nonstochastic generator $\lambda =\lambda _{+}^{-}$ as 
\begin{equation}
\boldsymbol{\lambda }\left( B\right) =\boldsymbol{L}^{*}\jmath (B)%
\boldsymbol{L}-\boldsymbol{K}^{*}B-B\boldsymbol{K}.  \label{0.6}
\end{equation}
The dilation of the stochastic differentials for CP processes over arbitrary 
$*$-algebras $\mathfrak{a}$ and $\mathcal{B}$, giving this structure for the
bounded generators over a von Neumann algebra $\mathcal{B}$ as a consequence
of the Christensen-Evans theorem \cite{CrE}, was obtained in \cite{15, Bge}.

Here we shall prove that the quantum stochastic extension (\ref{0.6}) of the
Lindblad's structure $\lambda \left( B\right) =L^{*}\jmath \left( B\right)
L-K^{*}B-BK$, can always be used for the construction and the dilation of
the CP flows also in the case of the unbounded maps $\lambda _\nu ^\mu $
over the algebra $\mathcal{B}=\mathcal{L}\left( \mathcal{H}\right) $. The
existence of minimal CP\ solution which has been recently constructed under
certain continuity conditions in \cite{Be97} proves that this structure is
also sufficient for the CP property of any solution to this stochastic
equation. The construction of the differential dilations and the CP
solutions of such quantum stochastic differential equations with the bounded
generators over the simple finite-dimensional It\^{o} algebra $\mathfrak{a}=%
\mathcal{Q}\left( \mathcal{E}\right) $ and the arbitrary $\mathcal{B}%
\subseteq \mathcal{L}\left( \mathcal{H}\right) $ was also recently discussed
in \cite{Be96, LiP}.

The nonstochastic case $\Lambda \left( t,a\right) =\alpha t\mathrm{I}$ is
described by the simplest, one-dimensional It\^{o} algebra $\mathfrak{a}=%
\mathbb{C}d $ with $l\left( a\right) =\alpha \in \mathbb{C}$ and the
nilpotent multiplication $\alpha ^{\star }\alpha =0$ corresponding to the
non-stochastic (Newton) calculus $\left( \mathrm{d}t\right) ^2=0$ in $%
\mathcal{E}=0$. The standard Wiener process $\mathrm{Q}=\Lambda
_{-}^{\bullet }+\Lambda _{\bullet }^{+}$ in Fock space is described by the
second order nilpotent algebra $\mathfrak{a}$ of pairs $a=\left( \alpha ,\xi
\right) $ with $d=\left( 1,0\right) $, $\xi \in \mathbb{C}$, represented by
the quadruples $a_{+}^{-}=\alpha ,\quad a_{\bullet }^{-}=\xi =a_{+}^{\bullet
},\quad a_{\bullet }^{\bullet }=0$ in $\mathcal{E}=\mathbb{C}$,
corresponding to $\Lambda \left( t,a\right) =\alpha t\mathrm{I}+\xi \mathrm{Q%
}\left( t\right) $. The unital $\star $-algebra $\mathbb{C}$ with the usual
multiplication $\zeta ^{\star }\zeta =\left| \zeta \right| ^2$ can be
embedded into the two-dimensional It\^{o} algebra $\mathfrak{a}$ of $%
a=\left( \alpha ,\zeta \right) $, $\alpha =l\left( a\right) $, $\zeta \in 
\mathbb{C}$ as $a_{\bullet }^{\bullet }=\zeta $, $a_{+}^{\bullet }=+i\zeta $%
, $a_{\bullet }^{-}=-i\zeta $, $a_{+}^{-}=\zeta $. It corresponds to $%
\Lambda \left( t,a\right) =\alpha t\mathrm{I}+\zeta \mathrm{P}\left(
t\right) $, where $\mathrm{P}=\Lambda _{\bullet }^{\bullet }+i\left( \Lambda
_{\bullet }^{+}-\Lambda _{-}^{\bullet }\right) $ is the representation of
the standard Poisson process, compensated by its mean value $t$. These two
commutative cases exhaust the possible types of two-dimensional It\^{o}
algebras. Thus, our results \cite{Bge, Be97} are also applicable to the
classical stochastic differentials of completely positive processes,
corresponding to the commutative It\^{o} algebras, which are always
decomposable into the Wiener, Poisson and Newton orthogonal components.

\section{Quantum sub-filtering dynamics}

The quantum filtering theory, which was outlined in \cite{1, 2} and
developed then since \cite{3}, provides the derivations for new types of
irreversible stochastic equations for quantum states, giving the dynamical
solution for the well-known quantum measurement problem. Some particular
types of such equations have been considered also in the phenomenological
theories of quantum permanent reduction \cite{4,5}, continuous measurement
collapse \cite{6,7}, spontaneous jumps \cite{8,9}, diffusions and
localizations \cite{10,11}. The main feature of such dynamics is that the
reduced irreversible evolution can be described in terms of a linear
dissipative stochastic wave equation, the solution to which is normalized
only in the mean square sense.

The simplest dynamics of this kind is described by the continuous filtering
wave propagators $V_t\left( \omega \right) $, defined on the space $\Omega $
of all Brownian trajectories as an adapted operator-valued stochastic
process in the system Hilbert space $\mathcal{H}$, satisfying the stochastic
diffusion equation 
\begin{equation}
\mathrm{d}V_t+KV_t\mathrm{d}t=LV_t\mathrm{dQ},\quad V_0=I  \label{1.3}
\end{equation}
in the It\^{o} sense. Here $\mathrm{Q}\left( t,\omega \right) $ is the
standard Wiener process, which is described by the independent increments $%
\mathrm{dQ}\left( t\right) =\mathrm{Q}\left( t+\mathrm{d}t\right) -\mathrm{Q}%
\left( t\right) $, having the zero mean values $\langle \mathrm{dQ}\rangle
=0 $ and the multiplication property $(\mathrm{dQ})^2=\mathrm{d}t$, $K$ is
an accretive operator, $K+K^{\dagger }\geq L^{*}L$, defined on a dense
domain $\mathcal{D\subseteq H}$, with $K^{\dagger }=K^{*}|\mathcal{D}$, and $%
L$ is a linear operator $\mathcal{D}\rightarrow \mathcal{H}$. This
stochastic wave equation with $K+K^{\dagger }=L^{*}L$ was first derived \cite%
{13} from a unitary cocycle evolution by a quantum filtering procedure. A
sufficient analyticity condition, under which it has the unique solution in
the form of stochastic multiple integral even in the case of unbounded $K$
and $L$ is given in \cite{Be97}. Using the It\^{o} formula 
\begin{equation}
\mathrm{d}\left( V_t^{*}V_t\right) =\mathrm{d}V_t^{*}V_t+V_t^{*}\mathrm{d}%
V_t+\mathrm{d}V_t^{*}\mathrm{d}V_t,  \label{1.2}
\end{equation}
and averaging $\left\langle \cdot \right\rangle $ over the trajectories of $%
\mathrm{Q}$, one obtains $\langle V_t^{*}V_t\rangle \leq I$ as a consequence
of $\mathrm{d}\langle V_t^{*}V_t\rangle \leq 0$. Note that the process $V_t$
is not necessarily unitary if the filtering condition $K^{\dagger }+K=L^{*}L$
holds, and even if $L^{\dagger }=-L$, it might be only isometric, $%
V_t^{*}V_t=I$, in the unbounded case.

Another type of the filtering wave propagator $V_t\left( \omega \right)
:\psi _0\in \mathcal{H}\mapsto \psi _t\left( \omega \right) $ in $\mathcal{H}
$ is given by the stochastic jump equation 
\begin{equation}
\mathrm{d}V_t+KV_t\mathrm{d}t=LV_t\mathrm{dP},\quad V_0=I.  \label{1.1}
\end{equation}
at the random time instants $\omega =\left\{ t_1,t_2,...\right\} $. Here $%
L=J-I$ is the jump operator, corresponding to the stationary discontinuous
evolutions $\psi _{t+}=J\psi $ at $t\in \omega $, and $\mathrm{P}\left(
t,\omega \right) =\mathrm{N}\left( t,\omega \right) -t$ is the standard
Poisson process, counting the number $\mathrm{N}\left( t,\omega \right)
=\left| \omega \cap [0,t)\right| $ compensated by its mean value $t$. It is
described as the process with independent increments $\mathrm{dP}\left(
t\right) =\mathrm{P}\left( t+\mathrm{d}t\right) -\mathrm{P}\left( t\right) $%
, having the values $\left\{ 0,1\right\} $ at $\mathrm{d}t\rightarrow 0$,
with zero mean $\langle \mathrm{dP}\rangle =0$, and the multiplication
property $\left( \mathrm{dP}\right) ^2=\mathrm{dP}+\mathrm{d}t$. This
stochastic wave equation was first derived in \cite{12} under the filtering
condition $L^{*}L=K+K^{\dagger }$ by the conditioning with respect to the
spontaneous reductions $J:\psi _t\mapsto \psi _{t+}$. An analyticity
condition under which it has the unique solution in the form of the multiple
stochastic integral even in the case of unbounded $K$ and $L$ is also given
in \cite{Be97}. Using the It\^{o} formula (\ref{1.2}) with $\mathrm{d}V_t^{*}%
\mathrm{d}V_t=V_t^{*}L^{*}LV_t(\mathrm{dP}\mathbf{+}\mathrm{d}t\mathrm{)}$,
one can obtain 
\begin{equation*}
\mathrm{d}\left( V_t^{*}V_t\right) =V_t^{*}\left( L^{*}L-K-K^{\dagger
}\right) V_t\mathrm{d}t+V_t^{*}\left( L^{\dagger }+L+L^{*}L\right) V_t%
\mathrm{dP}.
\end{equation*}
Averaging $\left\langle \cdot \right\rangle $ over the trajectories of $%
\mathrm{P}$, one can easily find that $\mathrm{d}\langle V_t^{*}V_t\rangle
\leq 0$ under the sub-filtering condition $L^{*}L\leq K+K^{\dagger }$. Such
evolution is not needed to be unitary, but in the filtering case it might be
isometric, $V_t^{*}V_t=I$ if the jumps are isometric, $J^{*}J=I$.

This proves in both cases that the stochastic wave function $\psi _t\left(
\omega \right) =V_t\left( \omega \right) \psi _0$ is not normalized for each 
$\omega $, but it is normalized in the mean square sense to the survival
probability $\langle ||\psi _t||^2\rangle \leq ||\psi _0||^2=1$, a
decreasing probability for a quantum unstable system not to be demolished
during its observation up to the time $t$. In the stable case $\left\langle
||\psi _t||^2\right\rangle =1$ the positive stochastic function $p_t\left(
\omega \right) =||\psi _t\left( \omega \right) ||^2$ is the probability
density of a diffusive $\widehat{\mathrm{Q}}$ or counting $\widehat{\mathrm{P%
}}$ output process up to the given $t$ with respect to the standard Wiener $%
\mathrm{Q}$ or Poisson $\mathrm{P}$ input processes correspondingly, in the
general case this is given by the conditional probability density $||\psi
_t\left( \omega \right) ||^2/\langle ||\psi _t||^2\rangle $.

Using the It\^{o} formula for $\rho _t\left( \omega \right) =V_t\left(
\omega \right) \rho _0V_t\left( \omega \right) ^{*}$, one can obtain the
stochastic equations 
\begin{equation}
\mathrm{d}\rho _t+\left( K\rho _t+\rho _tK^{*}-L\rho _tL^{*}\right) \mathrm{d%
}t=\left( L\rho _t+\rho _tL^{*}\right) \mathrm{dQ},  \label{1.5}
\end{equation}
\begin{equation}
\mathrm{d}\rho _t+\left( K\rho _t+\rho _tK^{*}-L\rho _tL^{*}\right) \mathrm{d%
}t=\left( J\rho _tJ^{*}-\rho _t\right) \mathrm{dP},  \label{1.4}
\end{equation}
describing the stochastic evolution $\Phi _t:\rho _0\mapsto \rho _t$ of an
initially normalized density operator $\rho _0\geq 0$, $\mathrm{tr}\rho _0=1$
as the stochastic density operator $\rho _t\left( \omega \right) =\Phi
_t\left( \omega ,\rho _0\right) $, normalized to the probability density $%
p_t\left( \omega \right) =\mathrm{tr}\rho _t\left( \omega \right) $. The
stochastic dynamical maps $\Phi _t\left( \rho \right) =V_t\rho V_t^{*}$ are
obviously positive but in general irreversible if $V_t\left( \omega \right) $
are not unitary, although they preserve the pure states in this particular
case.

Although the filtering equations (\ref{1.1}), (\ref{1.3}) look very
different, they can be unified in the form of quantum stochastic equation 
\begin{equation}
\mathrm{d}V_t+KV_t\mathrm{d}t+K^{-}V_t\mathrm{d}\Lambda _{-}=\left(
J-I\right) V_t\mathrm{d}\Lambda +L_{+}V_t\mathrm{d}\Lambda ^{+}  \label{1.6}
\end{equation}
where $\Lambda ^{+}\left( t\right) $ is the creation process, corresponding
to the annihilation $\Lambda _{-}\left( t\right) $ on the interval $[0,t)$,
and $\Lambda \left( t\right) $ is the number of quanta on this interval.
Indeed, the standard Poisson process $\mathrm{P}$ as well as the Wiener
process $\mathrm{Q}$ can be represented in $\mathfrak{F}$ by the linear
combinations \cite{16} 
\begin{equation}
\mathrm{P}\left( t\right) =\Lambda \left( t\right) +i\left( \Lambda
^{+}\left( t\right) -\Lambda _{-}\left( t\right) \right) ,\quad \mathrm{Q}%
\left( t\right) =\Lambda ^{+}\left( t\right) +\Lambda _{-}\left( t\right) ,
\label{1.9}
\end{equation}
so the equation (\ref{1.6}) corresponds to the stochastic diffusion equation
(\ref{1.3}) if $J=I$, $L_{+}=L=-K^{-}$, and it corresponds to the stochastic
jump equation (\ref{1.1}) if $J=I+L$, $L_{+}=iL=K^{-}$. These canonical
quantum stochastic processes, representing the quantum noise with respect to
the vacuum state $|0\rangle $ of the Fock space $\mathcal{F}$ over the
single-quantum Hilbert space $L^2\left( \mathbb{R}_{+}\right) $ of
square-integrable functions of $t\in [0,\infty )$, are formally given in 
\cite{14} by the integrals 
\begin{equation*}
\Lambda _{-}\left( t\right) =\int_0^t\Lambda _{-}^r\mathrm{d}r,\quad \Lambda
^{+}\left( t\right) =\int_0^t\Lambda _r^{+}\mathrm{d}r,\quad \Lambda \left(
t\right) =\int_0^t\Lambda _r^{+}\Lambda _{-}^r\mathrm{d}r,
\end{equation*}
where $\Lambda _{-}^r,\Lambda _r^{+}$ are the generalized quantum
one-dimensional fields in $\mathcal{F}$, satisfying the canonical
commutation relations 
\begin{equation*}
\left[ \Lambda _{-}^r,\Lambda _s^{+}\right] =\delta \left( s-r\right)
I,\quad \left[ \Lambda _{-}^r,\Lambda _{-}^s\right] =0=\left[ \Lambda
_r^{+},\Lambda _s^{+}\right] .
\end{equation*}
They can be defined by the independent increments with 
\begin{equation}
\langle 0|\mathrm{d}\Lambda _{-}|0\rangle =0,\quad \langle 0|\mathrm{d}%
\Lambda ^{+}|0\rangle =0,\quad \langle 0|\mathrm{d}\Lambda |0\rangle =0
\label{1.7}
\end{equation}
and the noncommutative multiplication table 
\begin{equation}
\mathrm{d}\Lambda \mathrm{d}\Lambda =\mathrm{d}\Lambda ,\quad \mathrm{d}%
\Lambda _{-}\mathrm{d}\Lambda =\mathrm{d}\Lambda _{-},\quad \mathrm{d}%
\Lambda \mathrm{d}\Lambda ^{+}=\mathrm{d}\Lambda ^{+},\quad \mathrm{d}%
\Lambda _{-}\mathrm{d}\Lambda ^{+}=\mathrm{d}tI  \label{1.8}
\end{equation}
with all other products being zero: $\mathrm{d}\Lambda \mathrm{d}\Lambda
_{-}=\mathrm{d}\Lambda ^{+}\mathrm{d}\Lambda =\mathrm{d}\Lambda ^{+}\mathrm{d%
}\Lambda _{-}=0$.

The corresponding quantum stochastic equation for the density operator $\rho
_t=V_t\rho _0V_t^{*}$ has the following form 
\begin{equation*}
\mathrm{d}\rho _t+\left( K\rho _t+\rho _tK^{*}-L\rho _tL^{*}\right) \mathrm{d%
}t=\left( J\rho _tJ^{*}-\rho _t\right) \mathrm{d}\Lambda
\end{equation*}
\begin{equation}
+\left( J\rho _tL^{-}-K^{-}\rho _t\right) \mathrm{d}\Lambda _{-}+\left(
L_{+}\rho _tJ^{*}-\rho _tK_{+}\right) \mathrm{d}\Lambda ^{+},  \label{1.10}
\end{equation}
where $L^{-}=L_{+}^{*},K_{+}^{*}=K^{-}$. The equation (\ref{1.10}),
coinciding with either (\ref{1.5}) or with (\ref{1.4}) in the particular
cases, is obtained from (\ref{1.6}) by using the It\^{o} formula (\ref{1.2})
with the multiplication table (\ref{1.8}). In the another particular case 
\begin{equation*}
J=S,\quad K^{-}=L^{-}S,\quad L_{+}=SK_{+},\quad S^{*}S=I,
\end{equation*}
it corresponds to the Hudson--Evans quantum stochastic flow \cite{EvH} if $%
S^{*}=S^{-1}$. Such evolution is isometric, and identity preserving, $%
V_tV_t^{*}=I$, i.e. unitary at least in the case of the bounded $K$ and $L$.

In the Heisenberg picture the stochastic dynamics is described by the dual
transformations $\phi _t\left( \omega \right) =\Phi _t^{\prime }\left(
\omega \right) $, such that for any density operator $\rho _0$ and for any
bounded observable $B$ on $\mathcal{H}$%
\begin{equation*}
\mathrm{tr}\left[ \Phi _t^{\prime }\left( \omega ,B\right) \rho _0\right] =%
\mathrm{tr}\left[ B\Phi _t\left( \omega ,\rho _0\right) \right] .
\end{equation*}
The linear stochastic maps $B\mapsto Y_t=\phi _t\left( B\right) $ are
obviously Hermitian in the sense that $Y_t^{*}=Y_t$ if $B^{*}=B$ and
completely positive, but in contrast to the usual Hamiltonian dynamics, they
are multiplicative, $\phi _t\left( B^{*}C\right) =\phi _t\left( B\right)
^{*}\phi _t\left( C\right) $ only in the case, corresponding to the HE flow,
even if they are not averaged with respect to $\omega $. Moreover, they are
usually not normalized, $R_t\left( \omega \right) :=\phi _t\left( \omega
,I\right) \neq I$, although the stochastic positive operators $%
R_t=V_t^{*}V_t $ under the filtering condition are usually normalized in the
mean, $\langle R_t\rangle =I$, and satisfy the martingale property $\epsilon
_t\left[ R_s\right] =R_t$ for all $s>t$, where $\epsilon _t$ is the
conditional expectation with respect to the history of the processes $%
\mathrm{P}$ or $\mathrm{Q}$ up to time $t$. The sub-filtering condition $%
K+K^{\dagger }\geq L^{-}L_{+}$ for the equation (\ref{1.6}) defines in both
cases the positive operator-valued stochastic process $R_t=\phi _t\left(
I\right) $ as a sub-martingale with $R_0=I$, or a martingale in the case $%
K+K^{\dagger }=L^{-}L_{+}$.

Although the filtering dynamics with unbounded coefficients of the
particular types has been studied elsewhere \cite{Ho96} by means of the
classical stochastic differential equations, the general structure of such
equations has not been discovered, and the general filtering CP flows have
not been constructed. In the next sections we define a multidimensional
analog of the quantum stochastic equation (\ref{1.10}), and will show that
the general structure of its generator indeed follows just from the property
of complete positivity of the dual stochastic maps $\phi _t=\Phi _t^{\prime
} $ for all $t>0$ and the normalization condition $\phi _t\left( I\right)
=R_t$ to a form-valued sub-martingale with respect to the natural filtration
of the quantum noise in the Fock space $\mathfrak{F}$ .

\section{Quantum completely positive flows}

Throughout the complex pre-Hilbert space $\mathcal{D}\subseteq \mathcal{H}$
is a Fr\'{e}chet (i.e. metrizable complete) space with respect to a stronger
topology, $\mathcal{E}\otimes \mathcal{D}$ denotes the projective tensor
product ($\pi $-product) with another such space $\mathcal{E}$, $\mathcal{D}%
^{\prime }\supseteq \mathcal{H}$ denotes the dual space of continuous
antilinear functionals $\eta ^{\prime }:\eta \in \mathcal{D}\mapsto \langle
\eta |\eta ^{\prime }\rangle $, with respect to the canonical pairing $%
\left\langle \eta |\eta ^{\prime }\right\rangle $ given by $\left\| \eta
\right\| ^2$ for $\eta ^{\prime }=\eta \in \mathcal{H}$, $\mathcal{B}\left( 
\mathcal{D}\right) $ denotes the linear space of all continuous sesquilinear
forms $\langle \eta |B\eta \rangle $ on $\mathcal{D}$, identified with the
continuous linear operators $B:\mathcal{D}\rightarrow \mathcal{D}^{\prime }$
(kernels), $B^{\dagger }\in \mathcal{B}\left( \mathcal{D}\right) $ is the
Hermit conjugated form (kernel) $\langle \eta |B^{\dagger }\eta \rangle
=\langle \eta |B\eta \rangle ^{*}$, and $\mathcal{L}\left( \mathcal{D}%
\right) \subseteq \mathcal{B}\left( \mathcal{D}\right) $ denotes the algebra
of all strongly continuous operators $B:\mathcal{D}\rightarrow \mathcal{D}$.
Any such space $\mathcal{D}$ can be considered as a projective limit with
respect to an increasing sequence of norms $\left\| \cdot \right\|
_p>\left\| \cdot \right\| $ on $\mathcal{D}$; for the definitions and
properties of this standard topological notions see for example \cite{Ob94}.
The spaces $\mathcal{D}^{\prime }$ and $\mathcal{B}\left( \mathcal{D}\right) 
$ will be equipped with w*- topologies induced by their preduals $\mathcal{D}
$ and $\mathcal{D}\otimes \mathcal{D}$, and coinciding with the weak
topology on each bounded subset with respect to a norm $\left\| \cdot
\right\| _p$. Any operator $A\in \mathcal{L}\left( \mathcal{D}\right) $ with 
$A^{\dagger }\in \mathcal{L}\left( \mathcal{D}\right) $ can be uniquely
extended to a weakly continuous operator onto $\mathcal{D}^{\prime }$ as $%
A^{\dagger *}$, denoted again as $A$, where $A^{*}$ is the dual operator $%
\mathcal{D}^{\prime }\rightarrow \mathcal{D}^{\prime }$, $\langle \eta
|A^{*}\eta ^{\prime }\rangle =\left\langle A\eta |\eta ^{\prime
}\right\rangle $, defining the involution $A\mapsto A^{*}$ for the
continuations $A:\mathcal{D}^{\prime }\rightarrow \mathcal{D}^{\prime }$. We
say that the operator $A$ commutes with a sesquilinear form, $BA=AB$ if $%
\left\langle \eta |BA\eta \right\rangle =\left\langle A^{\dagger }\eta
|B\eta \right\rangle $ for all $\eta \in \mathcal{D}$. The commutant $%
\mathcal{A}^c=\left\{ B\in \mathcal{B}\left( \mathcal{D}\right) :\left[ A,B%
\right] =0,\forall A\in \mathcal{A}\right\} $ of an operator $*$-algebra $%
\mathcal{A}\subseteq \mathcal{L}\left( \mathcal{D}\right) $ is weakly closed
in $\mathcal{B}\left( \mathcal{D}\right) $, so that the weak closure $%
\overline{\mathcal{B}}\subseteq \mathcal{B}\left( \mathcal{D}\right) $ of
any $\mathcal{B}\subseteq \mathcal{A}^c$ also commutes with $\mathcal{A}$.

Let us denote $\mathcal{B}=\mathcal{L}\left( \mathcal{H}\right) $ the
algebra of all bounded operators $B:\mathcal{H}\rightarrow \mathcal{H}$, $%
\left\| B\right\| <\infty $, $\overline{\mathcal{B}}=\mathcal{B}\left( 
\mathcal{D}\right) $ means the weak closure of $\mathcal{B}\subseteq 
\mathcal{B}\left( \mathcal{D}\right) $, and let $\left( \Omega ,\mathfrak{A}%
,P\right) $ be a probability space with a filtration $\left( \mathfrak{A}%
_t\right) _{t>0},$ $\mathfrak{A}_t\subseteq \mathfrak{A}$ of $\sigma $%
-algebras on $\Omega $. One can assume that the filtration $\mathfrak{A}%
_t\subseteq \mathfrak{A}_s,\forall t<s$ is generated by the pieces $%
x_{t]}=\left\{ r\mapsto x\left( r\right) :r\leq t\right\} $ of a stochastic
process $x\left( t,\omega \right) $ with independent increments $\mathrm{d}%
x\left( t\right) =x\left( t+\Delta \right) -x\left( t\right) $, and the
probability measure $P$ is invariant under the measurable representations $%
\omega \mapsto \omega _s\in \Omega $, $A_s^{-1}=\left\{ \omega :\omega _s\in
A\right\} \in \mathfrak{A}$, $\forall A\in \mathfrak{A}$ of the time shifts $%
t\mapsto t+s,s>0$ on $\Omega \ni \omega $, corresponding to the shifts of
the random increments 
\begin{equation*}
\mathrm{d}x\left( t,\omega _s\right) =\mathrm{d}x\left( t+s,\omega \right)
,\quad \forall \omega \in \Omega ,t\in \mathbb{R}_{+}.
\end{equation*}
The\emph{\ stochastic dynamics} over $\mathcal{B}$ with respect to the
process $x\left( t\right) $ is described by a cocycle flow $\phi =\left(
\phi _t\right) _{t>0}$ of linear completely positive \cite{Stn}
w*-continuous stochastic adapted maps $\phi _t\left( \omega \right) :%
\mathcal{B}\rightarrow \overline{\mathcal{B}}$, $\omega \in \Omega $ such
that the stochastic process $y_t\left( \omega \right) =\left\langle \eta
|\phi _t\left( \omega ,B\right) \eta \right\rangle $ is causally measurable
for each $\eta \in \mathcal{D}$, $B\in \mathcal{B}$ in the sense that $%
y_t^{-1}\left( B\right) \in \mathfrak{A}_t$, $\forall t>0$ and any Borel $%
B\subseteq \mathbb{C}$. The maps $\phi _t$ can be extended on the $\mathfrak{%
A}$-measurable functions $Y:\omega \mapsto Y\left( \omega \right) $ with
values $Y\left( \omega \right) \in \overline{\mathcal{B}}$ as the normal
maps $\phi _t\left[ Y\right] \left( \omega \right) =\overline{\phi }_t\left(
\omega ,Y\left( \omega _t\right) \right) $ for almost all $\omega \in \Omega 
$, where the linear maps $\overline{\phi }_t:\overline{\mathcal{B}}%
\rightarrow \overline{\mathcal{B}}$ are defined by the normal extensions of $%
\phi _t$ from the positive cone $\mathcal{B}_{+}$ onto $\overline{\mathcal{B}%
}_{+}$, so that the cocycle condition $\phi _r\left( \omega \right) \circ
\phi _s\left( \omega _r\right) =\phi _{r+s}\left( \omega \right) $, $\forall
r,s>0 $ reads as the semigroup condition $\phi _r\left[ \phi _s\left[ Y%
\right] \right] =\phi _{r+s}\left[ Y\right] $ of the extended maps. As it
was noted in the previous section, the maps $\phi _t\left( \omega \right) $
are not considered to be normalized to the identity, and can be even
unbounded, but in the case of filtering dynamics they are supposed to be
normalized, $\phi _t\left( \omega ,I\right) =R_t\left( \omega \right) $, to
an operator-valued martingale $R_t=\epsilon _t\left[ R_s\right] \geq 0$ with 
$R_0\left( \omega \right) =I$, or to a positive submartingale, $R_t\geq
\epsilon _t\left[ R_s\right] ,\forall s>t$ in the subfiltering case, where $%
\epsilon _t$ is the conditional expectation over $\omega $ with respect to $%
\mathfrak{A}_t$.

Now we give an algebraic generalization and a Fock space representation of
the filtering (or subfiltering) CP flows for a commutative It\^{o} algebra $%
\mathfrak{a}$, which was suggested in \cite{Be97} even in the noncommutative
case.

The role of the classical process $x\left( t\right) $ will play the quantum
stochastic process 
\begin{equation*}
X\left( t\right) =A\otimes I+I\otimes \Lambda \left( t,a\right) ,\quad A\in 
\mathcal{A},a\in \mathfrak{a}
\end{equation*}
parametrized by an Abelian $*$-subalgebra $\mathcal{A}\subset \mathcal{L}%
\left( \mathcal{D}\right) $ and a commutative It\^{o} algebra $\mathfrak{a}$%
. Here $\Lambda \left( t,a\right) $ is the process with independent
increment on a dense subspace $\mathfrak{F}\subset \Gamma \left( \mathfrak{E}%
\right) $ of the Fock space $\Gamma \left( \mathfrak{E}\right) $ over the
space $\mathfrak{E}=L_{\mathcal{E}}^2\left( \mathbb{R}_{+}\right) $ of all
square-norm integrable $\mathcal{E}$-valued functions on $\mathbb{R}_{+}$,
where $\mathcal{E}$ is a pre-Hilbert space of the representation $a\in 
\mathfrak{a}\mapsto \left( a_\nu ^\mu \right) _{\nu =+,\bullet }^{\mu
=-,\bullet }$ for the It\^{o} $\star $-algebra $\mathfrak{a}$. Every
commutative It\^{o} $*$-algebra is the sum $\mathfrak{a}=\mathbb{C}d+%
\mathfrak{a}_0+\mathfrak{a}_1$ of the Wiener and Poisson algebras $\mathbb{C}%
d+\mathfrak{a}_0$, $\mathbb{C}d+\mathfrak{a}_1$, such that each $a\in 
\mathfrak{a}$ has the unique decomposition $a=l\left( a\right) d+b+c$, where 
$d$ is the death of $\mathfrak{a}$, $b\in \mathfrak{a}_0$ is defined by the
conditions 
\begin{equation*}
b_{+}^{-}=l\left( b\right) =0,\qquad b_{\bullet }^{\bullet }=j\left(
b\right) =0,
\end{equation*}
and $c\in \mathfrak{a}_1$ is orthogonal to $b$: $bc=0=cb$, defined by the
condition $c_{\bullet }^{\bullet }=0\Rightarrow c_{+}^{\bullet
}=0=c_{\bullet }^{-}$. Thus the space $\mathcal{E}$ is decomposed into the
orthogonal sum $\mathcal{E}_0\oplus \mathcal{E}_1\oplus \mathcal{E}_{\bot }$
with $\mathcal{E}_0$ generated by $k\left( \mathfrak{a}_0\right) $, $%
\mathcal{E}_1$ generated by $k\left( \mathfrak{a}_1\right) $, and $\mathcal{E%
}_{\bot }$ is the orthogonal complement, which is zero if $\mathcal{E}$ is
the minimal space of the representation of $\mathfrak{a}$.

Assuming that $\mathcal{E}$ is a Fr\'{e}chet space, given by an increasing
sequence of Hilbertian norms $\left\| e^{\bullet }\right\| \left( \xi
\right) >\left\| e^{\bullet }\right\| $, $\xi \in \mathbb{N}$, we define $%
\mathfrak{F}$ as the projective limit $\cap _\xi \Gamma \left( \mathfrak{E}%
,\xi \right) $ of the Fock spaces $\Gamma \left( \mathfrak{E},\xi \right)
\subseteq \Gamma \left( \mathfrak{E}\right) $, generated by coherent vectors 
$f^{\otimes }$, with respect to the norms 
\begin{equation}
\left\| f^{\otimes }\right\| ^2\left( \xi \right) =\int_\Gamma \left\|
f^{\otimes }\left( \tau \right) \right\| ^2\left( \xi \right) \mathrm{d}\tau
:=\sum_{n=0}^\infty \frac 1{n!}\left( \int_0^\infty \left\| f^{\bullet
}\left( t\right) \right\| ^2\left( \xi \right) \mathrm{d}t\right)
^n=e^{\left\| f^{\bullet }\right\| ^2\left( \xi \right) }.  \label{2.0a}
\end{equation}
Here $f^{\otimes }\left( \tau \right) =\bigotimes_{t\in \tau }f^{\bullet
}\left( t\right) $ for each $f^{\bullet }\in \mathfrak{E}$ is represented by
tensor-functions on the space $\Gamma $ of all finite subsets $\tau =\left\{
t_1,...,t_n\right\} \subseteq \mathbb{R}_{+}$. Moreover, we shall assume
that the It\^{o} algebra $\mathfrak{a}$ is realized as a $\star $-subalgebra
of Hudson-Parthasarathy (HP) algebra $\mathcal{Q}\left( \mathcal{E}\right) $
of all quadruples $\boldsymbol{a}=\left( a_\nu ^\mu \right) _{\nu =+,\bullet
}^{\mu =-,\bullet }$ with $a_{\bullet }^{\bullet }\in \mathcal{L}\left( 
\mathcal{E}\right) $, strongly representing the $\star $-semigroup $1+%
\mathfrak{a}$ on the Fr\'{e}chet space $\mathcal{E}$ by projective
contractions $\delta _{\bullet }^{\bullet }+a_{\bullet }^{\bullet }\in 
\mathcal{L}\left( \mathcal{E}\right) $ in the sense that for each $\zeta \in 
\mathbb{N}$ there exists $\xi $ such that $\left\| e^{\bullet }+a_{\bullet
}^{\bullet }e^{\bullet }\right\| \left( \zeta \right) \leq \left\|
e^{\bullet }\right\| \left( \xi \right) $ for all $e^{\bullet }\in \mathcal{E%
}$. The following theorem proves that these are natural assumptions (which
are not restrictive in the simple Fock scale for a finite dimensional $%
\mathfrak{a}$.)

\begin{proposition}
The exponential operators $W\left( t,a\right) =:\exp \left[ \Lambda \left(
t,a\right) \right] :$ defined as the solutions to the quantum It\^{o}
equation 
\begin{equation}
\mathrm{d}W_t\left( g\right) =W_t\left( g\right) \mathrm{d}\Lambda \left(
t,g\left( t\right) \right) ,\quad W_0\left( g\right) =I,g\left( t\right) \in 
\mathfrak{a}  \label{2.1a}
\end{equation}
with $g\left( t\right) =a$, are strongly continuous, $W\left( t,a\right) \in 
\mathcal{L}\left( \mathfrak{F}\right) $, iff all $\widehat{a}_{\bullet
}^{\bullet }=\delta _{\bullet }^{\bullet }+a_{\bullet }^{\bullet }$ are
projective contructions on $\mathcal{E}$. They give an analytic
representation 
\begin{equation}
W\left( t,a\star a\right) =W\left( t,a\right) ^{*}W\left( t,a\right) ,\quad
W\left( t,0\right) =I,\quad W\left( t,d\right) =e^tI  \label{2.1g}
\end{equation}
of the unital $\star $-semigroup $1+\mathfrak{a}$ for the It\^{o} $\star $%
-algebra $\mathfrak{a}$ with respect to the $\star $-product $a\star
a=a+a^{\star }a+a^{\star }$.
\end{proposition}

\proof%
The solutions $W\left( t,a\right) $ are uniquely defined on the coherent
vectors as analytic functions 
\begin{equation}
W\left( t,a\right) f^{\otimes }\left( \tau \right) =\otimes _{r\in \tau
}^{r<t}\left( \widehat{a}_{\bullet }^{\bullet }f^{\bullet }\left( r\right)
+a_{+}^{\bullet }\right) \exp \left[ \int_0^t\left( a_{\bullet
}^{-}f^{\bullet }\left( r\right) +a_{+}^{-}\right) \mathrm{d}r\right]
\otimes _{r\in \tau }^{r\geq t}f^{\bullet }\left( r\right) ,  \label{2.1b}
\end{equation}
which obey the properties (\ref{2.1g}), see for example \cite{Bcs}. Thus the
span of coherent vectors is invariant, and it is also invariant under $%
W\left( t,a\right) ^{*}=W\left( t,a^{\star }\right) $. They can be extended
on $\mathfrak{F}$ by continuity which follows from the continuity of Wick
exponentials $\otimes \widehat{a}_{\bullet }^{\bullet }$ for the projective
contractions $\widehat{a}_{\bullet }^{\bullet }\in \mathcal{E}$ , and
boundedness of $a_{+}^{\bullet }\in \mathcal{E}$, $a_{\bullet }^{-}\in 
\mathcal{E}^{\prime }$. 
\endproof%

Let $\mathfrak{D}$ denote the Fr\'{e}chet space $\mathcal{D}\otimes 
\mathfrak{F}$, generated by $\psi =\eta \otimes f^{\otimes }$, $\eta \in 
\mathcal{D}$, $f^{\bullet }\in \mathfrak{E}$. Assuming for simplicity the
separability of the It\^{o} algebra in the sense $\mathcal{E}\subseteq \ell
^2$ such that $f^{\bullet }=\left( f^m\right) ^{m\in \mathbb{N}}$, one can
identify each $\psi ^{\prime }\in \mathfrak{D}^{\prime }$ with a sequence of 
$\mathcal{D}^{\prime }$-valued symmetric tensor-functions $\psi
_{m_1,...,m_n}^{\prime }\left( t_1,...t_n\right) $, $n=0,1,2,...$ . Let $%
\left( \mathfrak{D}_t\right) _{t>0}$ be the natural filtration and $\left( 
\mathfrak{D}_{[t}\right) _{t>0}$ be the backward filtration of the subspaces 
$\mathfrak{D}_t=\mathcal{D}\otimes \mathfrak{F}_t$, $\mathfrak{D}_{[t}=%
\mathcal{D}\otimes \mathfrak{F}_{[t}$ generated by $\eta \otimes f^{\otimes
} $ with $f^{\bullet }\in \mathfrak{E}_t$ and $f^{\bullet }\in \mathfrak{E}%
_{[t}$ respectively, where $\mathfrak{E}_t$ $=L_{\mathcal{E}}^2[0,t) $, $%
\mathfrak{E}_{[t}$ $=L_{\mathcal{E}}^2[t,\infty )$ are embedded into $%
\mathfrak{E}$. The spaces $\mathfrak{D}_t$, $\mathfrak{D}_{[t}$ of the
restrictions $E_t\psi =\psi |\Gamma _t$, $E_{[t}\psi =\psi |\Gamma _{[t}$
onto $\Gamma _t=\left\{ \tau _t=\tau \cap [0,t)\right\} $, $\Gamma
_{[t}=\left\{ \tau _{[t}=\tau \cap [t,\infty )\right\} $ are embedded into $%
\mathfrak{D}$ by the isometries $E_t^{\dagger }:\psi \mapsto \psi _t$, $%
E_{[t}^{\dagger }:\psi \mapsto \psi _{[t}$ as $\psi _t\left( \tau \right)
=\psi \left( \tau _t\right) \delta _\emptyset \left( \tau _{[t}\right) $, $%
\psi _{[t}\left( \tau \right) =\delta _\emptyset \left( \tau _t\right) \psi
\left( \tau _{[t}\right) $, where $\delta _\emptyset \left( \tau \right) =1$
if $\tau =\emptyset $, otherwise $\delta _\emptyset \left( \tau \right) =0$.
The projectors $E_t,E_{[t}$ onto $\mathfrak{D}_t,\mathfrak{D}^t$ are
extended onto $\mathfrak{D}^{\prime }$ as the adjoints to $E_t^{\dagger },$ $%
E_{[t}^{\dagger }$. The time shift on $\mathfrak{D}^{\prime }$ is defined by
the semigroup $\left( T^t\right) _{t>0} $ of adjoint operators $T^t=T_t^{*}$
to $T_t\psi \left( \tau \right) =\psi \left( \tau +t\right) $, where $\tau
+t=\left\{ t_1+t,...,t_n+t\right\} $, $\emptyset +t=\emptyset $, such that $%
T^t\psi \left( \tau \right) =\delta _\emptyset \left( \tau _t\right) \psi
\left( \tau _{[t}-t\right) $ are isometries for $\psi \in \mathfrak{D}$ onto 
$\mathfrak{D}_{[t}$. A family $\left( Z_t\right) _{t>0}$ of sesquilinear
forms $\left\langle \psi |Z_t\psi \right\rangle $ given by linear operators $%
Z_t:\mathfrak{D}\rightarrow \mathfrak{D}^{\prime }$ is called \emph{adapted}
(and $\left( Z^t\right) _{t>0}$ is called \emph{backward adapted}) if 
\begin{equation}
Z_t\left( \eta \otimes f^{\otimes }\right) =\psi ^{\prime }\otimes
E_{[t}f^{\otimes }\quad \left( Z^t\left( \eta \otimes f^{\otimes }\right)
=\psi ^{\prime }\otimes E_tf^{\otimes }\right) ,\quad \forall \eta \in 
\mathcal{D},f^{\bullet }\in \mathfrak{E},  \label{adp}
\end{equation}
where $\psi ^{\prime }\in \mathfrak{D}_t^{\prime }$ ($\mathfrak{D}%
_{[t}^{\prime }$) and $E_{[t}\;$($E_t$) are the projectors onto $\mathfrak{F}%
_{[t}$ ($\mathfrak{F}_t$) correspondingly.

The (vacuum) \emph{conditional expectation} on $\mathcal{B}\left( \mathfrak{D%
}\right) $ with respect to the past up to a time $t\in \mathbb{R}_{+}$ is
defined as a positive projector, $\epsilon _t\left( Z\right) \geq 0$, if $%
Z\geq 0$, $\epsilon _t=\epsilon _t\circ \epsilon _s,\forall s>t$, giving an
adapted sesquilinear form $Z_t=\epsilon _t\left( Z\right) $ in (\ref{adp})
for each $Z\in \mathcal{B}\left( \mathfrak{D}\right) $ by $\psi ^{\prime
}=E_tZE_t^{\dagger }\psi $, where $\psi =\eta \otimes E_tf^{\otimes }$. The
time shift $\left( \theta ^t\right) _{t>0}$ on $\mathcal{B}\left( \mathfrak{D%
}\right) $ is uniquely defined by the covariance condition $\theta ^t\left(
Z\right) T^t=T^tZ$ as a backward adapted family $Z^t=\theta ^t\left(
Z\right) ,t>0$ for each $Z\in \mathcal{B}\left( \mathfrak{D}\right) $. As in
the bounded case \cite{Mey} between the maps $\epsilon _t$ and $\theta ^t$
we have the relation $\theta ^r\circ \epsilon _s=\epsilon _{r+s}\circ \theta
^r$ which follows from the operator relation $T^rE_s=E_{r+s}T^r$. An adapted
family $\left( M_t\right) _{t>0}$ of positive $\left\langle \psi |M_t\psi
\right\rangle \geq 0,\forall \psi \in \mathfrak{D}$ Hermitian $M_t^{\dagger
}=M_t $ forms $M_t\in \mathcal{B}\left( \mathfrak{D}\right) $ is called 
\emph{martingale} (\emph{submartingale}) if $\epsilon _t\left( M_s\right)
=M_t$ ($\epsilon _t\left( M_s\right) \leq M_t$) for all $s\geq t\geq 0$.

Let $\mathfrak{B}$ denote the space of all $Y\in \mathcal{B}\left( \mathfrak{%
D}\right) $, commuting with all $X=\left\{ X\left( t\right) \right\} $ in
the sense 
\begin{equation*}
AY=YA,\quad \forall A\in \mathcal{A},\quad YW\left( t,a\right) =W\left(
t,a\right) Y,\quad \forall t>0,a\in \mathfrak{a},
\end{equation*}
where $A\left( \eta \otimes \varphi \right) =A\eta \otimes \varphi $, $%
W\left( \eta \otimes \varphi \right) =\eta \otimes W\varphi $, and the
unital $*$-algebra $\mathcal{B}\subseteq \mathcal{L}\left( \mathcal{H}%
\right) $ be weakly dense in the commutant $\mathcal{A}^c$ (we can take $%
\mathcal{B}=\mathcal{L}\left( \mathcal{H}\right) $ only if $\mathcal{A}=0$,
corresponding to $X\left( 0\right) =0$.) The quantum filtration $\left( 
\mathfrak{B}_t\right) _{t>0}$ is defined as the increasing family of
subspaces $\mathfrak{B}_t\subseteq \mathfrak{B}_s,t\leq s$ of the adapted
sesquilinear forms $Y_t\in \mathfrak{B}$. The covariant shifts $\theta
^t:Y\mapsto Y^t$ leave the space $\mathfrak{B}$ invariant, mapping it onto
the subspaces of backward adapted sesquilinear forms $Y^t=\theta ^t\left(
Y\right) $.

The \emph{quantum stochastic positive flow} over $\mathcal{B\ }$is described
by a one parameter family $\phi =\left( \phi _t\right) _{t>0}$ of linear
w*-continuous maps $\phi _t:$ $\mathcal{B}\rightarrow \mathfrak{B}$
satisfying

\begin{enumerate}
\item the causality condition $\phi _t\left( B\right) \subseteq \mathfrak{B}%
_t,\quad \forall B\in \mathcal{B},t\in \mathbb{R}_{+}$,

\item the complete positivity condition $\left[ \phi _t\left( B_{kl}\right) %
\right] \geq 0$ for each $t>0$ and for any positive definite matrix $\left[
B_{kl}\right] \geq 0$ with $B_{kl}\in \mathcal{B}$,

\item the cocycle condition $\phi _r\circ \phi _s^r=\phi _{r+s},\forall
t,s>0 $ with respect to the covariant shift $\phi _s^r=\theta ^r\circ \phi _s
$.
\end{enumerate}

Here the composition $\circ $ is understood as $\phi _r\left[ \phi _s\left(
B\right) \right] =\phi _{r+s}\left( B\right) $ in terms of the linear normal
extensions of $\phi _t\left[ B\otimes Z\right] =\overline{\phi }_t\left(
B\right) Z^t$ to the CP maps $\mathfrak{B}\rightarrow $ $\mathfrak{B}$,
forming a one-parameter semigroup, where $B\in \overline{\mathcal{B}}$, $%
\overline{\phi }_t:\overline{\mathcal{B}}\rightarrow \mathfrak{B}_t$ are the
normal extensions of $\phi _t$, $Z^t=\theta ^t\left( Z\right) $, $Z\in 
\mathcal{B}\left( \mathfrak{F}\right) $. These can be defined like in
classical case as $\phi _t\left[ Y\right] \left( \bar{f}^{\bullet
},f^{\bullet }\right) =\overline{\phi }_t\left( \bar{f}^{\bullet },Y\left( 
\bar{f}_t^{\bullet },f_t^{\bullet }\right) ,f^{\bullet }\right) $ with $%
f_t^{\bullet }\left( r\right) =f^{\bullet }\left( t+r\right) $ by the
coherent matrix elements $Y\left( \bar{f}^{\bullet },f^{\bullet }\right)
=F^{*}YF$ for $Y\in \mathfrak{B}$ given by the continuous operators $F:\eta
\mapsto \psi _f=\eta \otimes f^{\otimes }$ , $\eta \in \mathcal{D}$ for each 
$f^{\bullet }\in \mathfrak{E}_t$ with the adjoints $F^{*}\psi ^{\prime
}=\int_{\tau <t}f^{\otimes }\left( \tau \right) ^{*}\psi ^{\prime }\left(
\tau \right) \mathrm{d}\tau $ for $\psi ^{\prime }\in \mathfrak{D}^{\prime }$%
.

The flow is called \emph{(sub)-filtering}, if $R_t=\phi _t\left( I\right) $
is a (sub)-martingale with $R_0=I$, and is called contractive, if $I\geq
R_t\geq R_s$ for all $0\leq t\leq s\in \mathbb{R}_{+}$.

\begin{proposition}
The complete positivity for adapted linear maps $\phi _t:\mathcal{B}%
\rightarrow \mathcal{B}\left( \mathfrak{D}\right) $ can be written as 
\begin{equation}
\sum_{f,h\in \mathfrak{E}_t}\sum_{B,C\in \mathcal{B}}\langle \xi _B^f|\phi
_t\left( \bar{f}^{\bullet },B^{*}C,h^{\bullet }\right) \xi _C^h\rangle
:=\langle \eta ^k|\phi _t\left( \bar{f}_k^{\bullet },B_k^{*}B_l,h_l^{\bullet
}\right) \eta ^l\rangle \geq 0,\quad \forall t>0  \label{cpd}
\end{equation}
(the usual summation rule over repeated cross-level indices is understood),
where $\xi _B^f=\eta ^k$ if $f^{\bullet }=f_k^{\bullet }$ and $B=B_k$ with $%
f_k^{\bullet }\in \mathfrak{E}_t,B_k\in \mathcal{B}$, $k=1,2,...,$ otherwise 
$\xi _B^f=0$, and $\phi _t\left( B,f^{\bullet }\right) =\phi _t\left(
B\right) F$, $\phi _t\left( \bar{f}^{\bullet },B\right) =F^{*}\phi _t\left(
B\right) $.
\end{proposition}

\proof%
By definition the map $\phi $ into the sesquilinear forms is completely
positive on $\mathcal{B}$ if $\left\langle \psi ^k|\phi \left( B_{kl}\right)
\psi ^l\right\rangle \geq 0$ whenever $\left\langle \eta ^k|B_{kl}\eta
^l\right\rangle \geq 0$, where $\eta ^k,\psi ^k$ are arbitrary finite
sequences. Approximating from below the latter positive forms by sums of the
forms $\sum_{kl}\left\langle \eta ^k|B_{ik}^{*}B_{il}\eta ^l\right\rangle
\geq 0$, the complete positivity can be tested only for the forms $%
\sum_{kl}\left\langle \eta ^k|B_k^{*}B_l\eta ^l\right\rangle \geq 0$ due to
the additivity $\phi \left( \sum_iB_{ik}^{*}B_{il}\right) =\sum_i\phi \left(
B_{ik}^{*}B_{il}\right) $. If $\phi _t$ is adapted, this can be written as 
\begin{equation*}
\sum_{B,C\in \mathcal{B}}\left\langle \chi _B|\phi \left( B^{*}C\right) \chi
_C\right\rangle =\left\langle \psi ^k|\phi \left( B_k^{*}B_l\right) \psi
^l\right\rangle :=\sum_{k,l}\left\langle \psi ^k|\phi \left(
B_k^{*}B_l\right) \psi ^l\right\rangle \geq 0,
\end{equation*}
where $\chi _B=\psi ^k\in \mathfrak{D}_t$ if $B=B_k\in \mathcal{B}$,
otherwise $\chi _B=0$. Because any $\psi \in \mathfrak{D}_t$ can be
approximated by a $\mathcal{D}$-span $\sum_f\eta ^f\otimes f^{\otimes }$of
coherent vectors over $f_k^{\bullet }\in \mathfrak{E}_t$, it is sufficient
to define the CP property only for such spans as 
\begin{equation*}
0\leq \sum_{f,h}\sum_{B,C}\left\langle \xi _B^f\otimes f^{\otimes }|\phi
\left( B^{*}C\right) \left( \xi _C^h\otimes h^{\otimes }\right)
\right\rangle =\sum_{f,h}\sum_{B,C}\left\langle \xi _B^f|\phi \left( \bar{f}%
^{\bullet },B^{*}C,h^{\bullet }\right) \xi _C^h\right\rangle .
\end{equation*}
\endproof%

Note that the subfiltering (filtering) flows can be considered as a quantum
stochastic CP dilations of the quantum sub-Markov (Markov) semigroups $%
\theta =\left( \theta _t\right) _{t>0}$, $\theta _r\circ \theta _s=\theta
_{r+s}$ in the sense $\theta _t=\epsilon \circ \phi _t$, where $\epsilon
\left( Y\right) \eta =EY\psi _0$, $E\psi ^{\prime }=\psi ^{\prime }\left(
\emptyset \right) ,\forall \psi ^{\prime }\in \mathfrak{D}^{\prime }$, with $%
\theta _s\left( I\right) \leq \theta _t\left( I\right) \leq I$ ($\theta
_t\left( I\right) =I$), $\forall t\leq s$. The contraction $C_t=\theta
_t\left( I\right) $ with $C_0=I$ defines the probability $\left\langle \eta
|C_t\eta \right\rangle \leq 1$, $\forall \eta \in \mathcal{H},\left\| \eta
\right\| =1$ for an unstable system not to be demolished by a time $t\in 
\mathbb{R}_{+}$, and the conditional expectations $\left\langle \eta
|AC_t\eta \right\rangle /\left\langle \eta |C_t\eta \right\rangle $ of the
initial nondemolition observables $A\in \mathcal{A}$ in any state $\eta \in 
\mathcal{D}$, and thus in any initial state $\psi _0\in \mathcal{\eta }%
\otimes \delta _\emptyset $. The following theorem shows that the
submartingale (or the contraction) $R_t=\phi _t\left( I\right) $ is the
density operator with respect to $\psi _0=\eta \otimes \delta _\emptyset $, $%
\eta \in \mathcal{H}$ (or with respect to any $\psi \in \mathcal{H}\otimes 
\mathfrak{F}$) also for the conditional state of the restricted
nondemolition process $X_{t]}=\left\{ r\mapsto X\left( r\right) :r\leq
t\right\} $.

\begin{theorem}
Let $t\mapsto R_t\in \mathfrak{B}_t$ be a positive (sub)-martingale and $%
\left( \mathfrak{g}_t\right) _{t>0}$ be the increasing family of $\star $%
-semigroups $\mathfrak{g}_t$ of step functions $g:\mathbb{R}_{+}\rightarrow 
\mathfrak{a}$, $g\left( s\right) =0$, $\forall s\geq t$ under the $\star $%
-product 
\begin{equation}
\left( g_k\star g_l\right) \left( t\right) =g_l\left( t\right) +g_k\left(
t\right) ^{\star }g_l\left( t\right) +g_k\left( t\right) ^{\star }
\label{2.1c}
\end{equation}
of $g_k^{\star }=g_k\star 0$ and $g_l=0\star g_l$. The generating function $%
\vartheta _t\left( g\right) =\epsilon \left[ R_tW_t\left( g\right) \right] $
of the output state for the process $\Lambda \left( t\right) $, defined for
any $g\in \mathfrak{g}_t$ and each $t>0$ as 
\begin{equation}
\left\langle \eta |\vartheta _t\left( g\right) \eta \right\rangle
=\left\langle \psi _0|R_tW_t\left( g\right) \psi _0\right\rangle ,\quad \psi
_0=\eta \otimes \delta _\emptyset ,  \label{sgf}
\end{equation}
is $\mathcal{B}^c$-valued, positive, $\vartheta _t\geq 0$ in the sense of
positive definiteness of the kernel 
\begin{equation}
\left\langle \eta ^k|\vartheta _t\left( g_k\star g_l\right) \eta
^l\right\rangle \geq 0,\quad \forall g_k\in \mathfrak{g}_t;\eta ^k\in 
\mathcal{D},  \label{2.1d}
\end{equation}
and $\vartheta _t\geq \vartheta _s|\mathfrak{g}_t$ in this sense for any $%
s\geq t $. If $R_0=I$, then $\vartheta _0\left( 0\right) =I\geq \vartheta
_t\left( 0\right) $ , and if $R_t$ is a martingale, then $\vartheta
_t=\vartheta _s|\mathfrak{g}_t$ for any $s\geq t$, and $\vartheta _t\left(
0\right) =I$ for all $t\in \mathbb{R}_{+}$. Any family $\vartheta =\left(
\vartheta _t\right) _{t\geq 0}$ of positive-definite functions $\vartheta _t:%
\mathfrak{g}_t\rightarrow \mathcal{B}^c$, satisfying the above consistency
and normalization properties, is the state generating function of the form (%
\ref{sgf}) iff it is absolutely continuous in the following sense 
\begin{equation}
\lim_{n\rightarrow \infty }\sum_{g\in \mathfrak{g}_t}\eta _n^g\otimes
g_{+}^{\otimes }=0\Rightarrow \lim_{n\rightarrow \infty }\sum_{g,h\in 
\mathfrak{g}_t}\left\langle \eta _n^g|\vartheta _t\left( g\star h\right)
\eta _n^h\right\rangle =0,  \label{2.1e}
\end{equation}
where $g_{+}^{\otimes }\left( \tau \right) =\otimes _{t\in \tau
}g_{+}^{\bullet }\left( t\right) $ and $\eta _n^g=0$ for almost all $g$
(i.e. except for a finite number of $g\in \mathfrak{g}_t$).
\end{theorem}

The proof is given in \cite{Be97} even for the general (noncommutative)
algebras $\mathcal{A}$ and $\mathfrak{a}$.

\section{Generators of quantum CP dynamics}

The quantum stochastically differentiable positive flow $\phi $ is defined
as a weakly continuous function $t\mapsto \phi _t$ with CP values $\phi _t:%
\mathcal{B}\rightarrow \mathfrak{B}_t$, $\phi _0\left( B\right) =B\otimes
I,\forall B\in \mathcal{B}$ such that for any product-vector $\psi _f=\eta
\otimes f^{\otimes }$ given by $\eta \in \mathcal{D}$ and $f^{\bullet }\in 
\mathfrak{E}$ 
\begin{equation}
\frac{\mathrm{d}}{\mathrm{d}t}\left\langle \psi _f|\phi _t\left( B\right)
\psi _f\right\rangle =\left\langle \psi _f|\phi _t\left( \lambda \left( \bar{%
f}^{\bullet }\left( t\right) ,B,f^{\bullet }\left( t\right) \right) \right)
\psi _f\right\rangle ,\qquad B\in \mathcal{B},  \label{2.1}
\end{equation}
where $\lambda \left( \bar{e}^{\bullet },B,e^{\bullet }\right) =$ $\lambda
\left( B\right) +e_{\bullet }\lambda ^{\bullet }\left( B\right) +\lambda
_{\bullet }\left( B\right) e^{\bullet }+e_{\bullet }\lambda _{\bullet
}^{\bullet }\left( B\right) e^{\bullet }$, $e_{\bullet }=\bar{e}^{\bullet }$
is the linear form on $\mathcal{E}$ with $e_{\bullet }^{*}=e^{\bullet }\in 
\mathcal{E}$ and $\left\langle \psi _f|\phi _0\left( B\right) \psi
_f\right\rangle =\left\langle \eta |B\eta \right\rangle \exp \left\|
f^{\bullet }\right\| ^2$. The generator $\lambda \left( B\right) =\lambda
\left( 0,B,0\right) $ of the quantum dynamical semigroup $\theta _t=\epsilon
\circ \phi _t\,$ is a linear w*-continuous map $B\mapsto \lambda \left(
B\right) \in \mathcal{A}^c$, $\lambda ^{\bullet }=\lambda _{\bullet
}^{\dagger }$ is a linear w*-continuous map given by the Hermitian adjoint
values $\lambda _{\bullet }\left( B^{*}\right) =\lambda ^{\bullet }\left(
B\right) ^{\dagger }$ in the continuous operators $\mathcal{E}\rightarrow 
\mathcal{A}^c$, and $\lambda _{\bullet }^{\bullet }:\mathcal{B}\rightarrow 
\mathcal{B}\left( \mathcal{D}\otimes \mathcal{E}\right) $is a w*-continuous
map with the values $\lambda _{\bullet }^{\bullet }\left( B\right) $ given
by continuous operators $\mathcal{E}\otimes \mathcal{E}\rightarrow \mathcal{A%
}^c$. The differential evolution equation (\ref{2.1}) for the coherent
vector matrix elements $\left\langle \psi _f|\phi _t\left( B\right) \psi
_f\right\rangle $ corresponds to the It\^{o} form \cite{16} of the quantum
stochastic equation 
\begin{equation}
\mathrm{d}\phi _t\left( B\right) =\phi _t\circ \lambda _\nu ^\mu \left(
B\right) \mathrm{d}\Lambda _\mu ^\nu :=\sum_{\mu ,\nu }\phi _t\left( \lambda
_\nu ^\mu \left( B\right) \right) \mathrm{d}\Lambda _\mu ^\nu ,\qquad \text{ 
}B\in \mathcal{B}\qquad  \label{2.2}
\end{equation}
with the initial condition $\phi _0\left( B\right) =B$, for all $B\in 
\mathcal{B}$. Here $\lambda _\nu ^\mu $ are the flow generators $\lambda
_{+}^{-}=\lambda $, $\lambda _{+}^{\bullet }=\lambda ^{\bullet }$, $\lambda
_{\bullet }^{-}=\lambda _{\bullet }$, $\lambda _{\bullet }^{\bullet }$,
called the structural maps, and the summation is taken over the indices $\mu
=-,\bullet $, $\nu =+,\bullet $ of the standard quantum stochastic
integrators $\Lambda _\mu ^\nu $. For simplicity we shall assume that the
pre-Hilbert Fr\'{e}chet space $\mathcal{E}$ is separable, $\mathcal{E}%
\subseteq $ $\ell ^2$. Then the index $\bullet $ can take any value in $%
\left\{ 1,2,...\right\} $ and $\Lambda _\mu ^\nu \left( t\right) $ are
indexed with $\mu \in \left\{ -,1,2,...\right\} $, $\nu \in \left\{
+,1,2,...\right\} $ as the standard time $\Lambda _{-}^{+}\left( t\right) =t%
\mathrm{I}$, annihilation $\Lambda _{-}^m\left( t\right) $, creation $%
\Lambda _n^{+}\left( t\right) $ and exchange-number $\Lambda _n^m\left(
t\right) $ operator integrators with $m,n\in \mathbb{N}$. The infinitesimal
increments $\mathrm{d}\Lambda _\nu ^\mu \left( t\right) =\Lambda _\nu ^{t\mu
}\left( \mathrm{d}t\right) $ are formally defined by the HP multiplication
table \cite{16} and the $\star $ -property \cite{3}, 
\begin{equation}
\mathrm{d}\Lambda _\mu ^\alpha \mathrm{d}\Lambda _\beta ^\nu =\delta _\beta
^\alpha \mathrm{d}\Lambda _\mu ^\nu ,\qquad \text{ }\Lambda ^{\star
}=\Lambda ,\qquad  \label{2.3}
\end{equation}
where $\delta _\beta ^\alpha $ is the usual Kronecker delta restricted to
the indices $\alpha \in \left\{ -,1,2,...\right\} $, $\beta \in \left\{
+,1,2,...\right\} $ and $\Lambda _{-\nu }^{\star \mu }=\Lambda _{-\mu }^{\nu
*}$ with respect to the reflection $-(-)=+,$ $-(+)=-$ of the indices $\left(
-,+\right) $ only.

The linear equation (\ref{2.2}) of a particular type, (quantum Langevin
equation) with bounded finite-dimensional structural maps $\lambda _\nu ^\mu 
$ was introduced by Evans and Hudson \cite{EvH} in order to describe the $*$%
-homomorphic quantum stochastic evolutions. The constructed quantum
stochastic $*$-homomorphic flow (EH-flow) is identity preserving and is
obviously completely positive, but it is hard to prove these algebraic
properties for the unbounded case. However the typical quantum filtering
dynamics is not homomorphic or identity preserving, but it is completely
positive and in the most interesting cases is described by unbounded
generators $\lambda _\nu ^\mu $. In the general content the equation (\ref%
{2.2}) was studied in \cite{20}, and the correspondent quantum stochastic,
not necessarily homomorphic and normalized flow was constructed even for the
infinitely-dimensional non-adapted case under the natural integrability
condition for the chronological products of the generators $\lambda _\nu
^\mu $ in the norm scale (\ref{2.0a}). The EH flows with unbounded $\lambda
_\nu ^\mu $, satisfying certain analyticity conditions, have been recently
constructed in strong sense by Fagnola-Sinha in \cite{18} for the
non-Hilbert class $L^\infty $ of test functions $f^{\bullet }$. Here we will
formulate the necessary differential conditions which follow from the
complete positivity, causality, and martingale properties of the filtering
flows, and which are sufficient for the construction of the quantum
stochastic flows obeying these properties in the case of the bounded $%
\lambda _\nu ^\mu $. As we showed in \cite{15}, the found properties are
sufficient to define the general structure of the bounded generators, and
this structure will help us in construction of the minimal completely
positive weak solutions for the quantum filtering equations also with
unbounded $\lambda _\nu ^\mu $.

Obviously the linear w*-continuous generators $\lambda _\nu ^\mu :\mathcal{B}%
\rightarrow \mathcal{A}^c$ for CP flows $\phi _t^{*}=\phi _t$, where $\phi
_t^{*}\left( B\right) =\phi _t\left( B^{*}\right) ^{\dagger }$, must satisfy
the $\star $ -property $\lambda ^{\star }=\lambda $, where $\lambda _{-\mu
}^{\star \nu }=\lambda _{-\nu }^{\mu *}$, $\lambda _\nu ^{\mu *}\left(
B\right) =\lambda _\nu ^\mu \left( B^{*}\right) ^{*}$ and are independent of 
$t$, corresponding to cocycle property $\phi _s\circ \phi _r^s=\phi _{s+r}$,
where $\phi _t^s$ is the solution to (\ref{2.2}) with $\Lambda _\nu ^\mu
\left( t\right) $ replaced by $\Lambda _\nu ^{s\mu }\left( t\right) $, and $%
\lambda _{+}^{-}\left( I\right) =0$ if $\phi $ is a filtering flow, $\phi
_t\left( I\right) =I$, as it is in the multiplicative case \cite{EvH}. We
shall assume that $\boldsymbol{\lambda }=\left( \lambda _\nu ^\mu \right)
_{\nu =+,\bullet }^{\mu =-,\bullet }$ for each $B^{*}=B$ defines a
continuous Hermitian form $\boldsymbol{b}=\boldsymbol{\lambda }\left(
B\right) $ on the Fr\'{e}chet space $\mathcal{D}\oplus \mathcal{D}_{\bullet
} $, 
\begin{equation*}
\left\langle \boldsymbol{\eta }\right| \boldsymbol{b}\left. \boldsymbol{\eta 
}\right\rangle =\sum_{m,n}\left\langle \eta ^m|b_n^m\eta ^n\right\rangle
+\sum_m\left\langle \eta ^m|b_{+}^m\eta \right\rangle +\sum_n\left\langle
\eta |b_n^{-}\eta ^n\right\rangle +\left\langle \eta |b_{+}^{-}\eta
\right\rangle ,
\end{equation*}
where $\eta \in \mathcal{D}$, $\eta ^{\bullet }=\left( \eta ^m\right) ^{m\in 
\mathbb{N}}\in \mathcal{D}_{\bullet }=\mathcal{D}\otimes \mathcal{E}$. We
say that an It\^{o} algebra $\mathfrak{a}$ , represented on $\mathcal{E}$,
commutes in HP sense with a $\boldsymbol{b}$, given by the form-generator $%
\boldsymbol{\lambda }$ if $\left( I\otimes a_{\bullet }^\mu \right) b_\nu
^{\bullet }=b_{\bullet }^\mu \left( I\otimes a_\nu ^{\bullet }\right) $ (For
simplicity the ampliation $I\otimes a_\nu ^\mu $ will be written again as $%
a_\nu ^\mu $.) Note that if we define the matrix elements $a_\nu ^\mu $, $%
b_\nu ^\mu $ also for $\mu =+$ and $\nu =-$, by the extension 
\begin{equation*}
a_\nu ^{+}=0=a_{-}^\mu ,\qquad \lambda _\nu ^{+}\left( B\right) =0=\lambda
_{-}^\mu \left( B\right) ,\quad \forall a\in \mathfrak{a},B\in \mathcal{B},
\end{equation*}
the HP product (\ref{0.4}) of $\boldsymbol{a}$ and $\boldsymbol{b}$ can be
written in terms of the usual matrix product $\mathbf{ab}=\left[ a_\lambda
^\mu b_\nu ^\lambda \right] $ of the extended quadratic matrices $\mathbf{a}=%
\left[ a_\nu ^\mu \right] _{\nu =-,\bullet ,+}^{\mu =-,\bullet ,+}$ and $%
\mathbf{b=}\boldsymbol{b}\mathbf{g}$, where $\mathbf{g}=\left[ \delta _{-\nu
}^\mu \right] $. Then one can extend the summation in (\ref{2.2}) so it is
also over $\mu =+$, and $\nu =-$, such that $b_\nu ^\mu \mathrm{d}\Lambda
_\mu ^\nu $ is written as the trace $\mathbf{b\cdot }\mathrm{d}\boldsymbol{%
\Lambda }$ over all $\mu ,\nu $. By such an extension the multiplication
table for $\mathrm{d}\Lambda \left( a\right) =\mathbf{a\cdot }\mathrm{d}%
\boldsymbol{\Lambda }$ , $\mathrm{d}\Lambda \left( b\right) =\mathbf{b\cdot }%
\mathrm{d}\boldsymbol{\Lambda }$ can be represented as $\mathrm{d}\Lambda
\left( a\right) \mathrm{d}\Lambda \left( b\right) =\mathbf{ab\cdot }\mathrm{d%
}\boldsymbol{\Lambda }$, and the involution $\mathbf{b\mapsto b}^{\star }$,
defining $\mathrm{d}\Lambda \left( b\right) ^{\dagger }=\mathbf{b}^{\star }%
\mathbf{\cdot }\mathrm{d}\boldsymbol{\Lambda }$, can be obtained by the
pseudo-Hermitian conjugation $b_\alpha ^{\star \nu }=g_{\alpha \mu }b_\beta
^{\mu *}g^{\beta \nu }$ respectively to the indefinite Minkowski metric
tensor $\mathbf{g}=\left[ g_{\mu \nu }\right] $ and its inverse $\mathbf{g}%
^{-1}=\left[ g^{\mu \nu }\right] $, given by $g^{\mu \nu }=\delta _{-\nu
}^\mu I=g_{\mu \nu }$.

Now let us find the differential form of the normalization and causality
conditions with respect to the quantum stationary process, with independent
increments $\mathrm{d}X\left( t\right) =X\left( t+\Delta \right) -X\left(
s\right) $ generated by an It\^{o} algebra $\mathfrak{a}$ on the separable
space $\mathcal{E}$.

\begin{proposition}
Let $\phi $ be a flow, satisfying the quantum stochastic equation (\ref{2.2}%
), and $\left[ W_t\left( g\right) ,\phi _t\left( B\right) \right] =0$ for
all $g\in \mathfrak{g},B\in \mathcal{B}$. Then the coefficients $b_\nu ^\mu
=\lambda _\nu ^\mu \left( B\right) $, $\mu =-,\bullet $, $\nu =+,\bullet $,
where $\bullet =1,2,...$, written in the matrix form $\boldsymbol{b}=\left(
b_\nu ^\mu \right) _{\nu =+,\bullet }^{\mu =-,\bullet }$, commute in the
sense of the HP product with $\boldsymbol{a}=\left( a_\nu ^\mu \right) _{\nu
=+,\bullet }^{\mu =-,\bullet }$ for all $a\in \mathfrak{a}$ and $B\in 
\mathcal{B} $: 
\begin{equation}
\left[ \boldsymbol{a},\boldsymbol{b}\right] :=\left( a_{\bullet }^\mu b_\nu
^{\bullet }-b_{\bullet }^\mu a_\nu ^{\bullet }\right) _{\nu =+,\bullet
}^{\mu =-,\bullet }=0.  \label{2.4}
\end{equation}
\end{proposition}

\proof%
Since $\epsilon _t\left( \phi _s\left( I\right) -\phi _t\left( I\right)
\right) $ is a negative Hermitian form, 
\begin{equation*}
\epsilon _t\left( \mathrm{d}\phi _t\left( I\right) \right) =\epsilon
_t\left( \phi _t\left( \lambda _\nu ^\mu \left( I\right) \right) \mathrm{d}%
\Lambda _\mu ^\nu \right) =\phi _t\left( \lambda _{+}^{-}\left( I\right)
\right) \mathrm{d}t\leq 0.
\end{equation*}
Since $Y_t=\phi _t\left( B\right) $ commutes with $W_t\left( g\right) $ for
all $B$ and $g\left( t\right) =a$, we have by virtue of quantum It\^{o}'s
formula 
\begin{equation*}
\mathrm{d}\left[ Y_t,W_t\right] =\left[ \mathrm{d}Y_t,W_t\right] +\left[ Y_t,%
\mathrm{d}W_t\right] +\left[ \mathrm{d}Y_t,\mathrm{d}W_t\right] =0.
\end{equation*}
The equations (\ref{2.1a}), (\ref{2.2}) and commutativity of $a_\nu ^\mu $
with $Y_t$ and $W_t$ imply 
\begin{eqnarray*}
&&\left( \left[ \phi _t\left( b_\nu ^\mu \right) ,W_t\right] +\left[
Y_t,a_\nu ^\mu W_t\right] +\phi _t\left( b_{\bullet }^\mu \right) a_\nu
^{\bullet }W_t-a_{\bullet }^\mu W_t\phi _t\left( b_\nu ^{\bullet }\right)
\right) \mathrm{d}\Lambda _\mu ^\nu \\
&=&W_t\left( \phi _t\left( b_{\bullet }^\mu \right) a_\nu ^{\bullet
}-a_{\bullet }^\mu \phi _t\left( b_\nu ^{\bullet }\right) \right) \mathrm{d}%
\Lambda _\nu ^\mu =W_t\phi _t\left( b_{\bullet }^\mu a_\nu ^{\bullet
}-a_{\bullet }^\mu b_\nu ^{\bullet }\right) \mathrm{d}\Lambda _{\mu .}^\nu
=0.
\end{eqnarray*}
Thus $\boldsymbol{a}\bullet \boldsymbol{b}=\boldsymbol{b}\bullet \boldsymbol{%
a}$ by the argument \cite{Par} of independence of the integrators $\mathrm{d}%
\Lambda _\mu ^\nu $. 
\endproof%

In order to formulate the CP differential condition we need the notion of 
\emph{quantum stochastic germ} for the CP flow $\phi $ at $t=0$. It was
defined in \cite{20, Bge}, for a quantum stochastic differential (\ref{2.2})
with $\phi _0\left( B\right) =B,\forall B\in \mathcal{B}$ as $\gamma _\nu
^\mu =\lambda _\nu ^\mu +\imath _\nu ^\mu $, where $\lambda _\nu ^\mu $ are
the structural maps $B\mapsto \lambda _\nu ^\mu \left( B\right) $ given by
the generators of the quantum It\^{o} equation (\ref{2.2}) and $\imath _\nu
^\mu :B\mapsto B\delta _\nu ^\mu $ is the ampliation of $\mathcal{B}$. Let
us prove that the germ-maps $\gamma _\nu ^\mu $ of a CP flow $\phi $ must be
conditionally completely positive (CCP) in a degenerated sense as it was
found for the finite-dimensional bounded case in \cite{15, Be96}. Another,
equivalent, but not so explicit characterization was suggested for this
particular case in \cite{LiP}.

\begin{theorem}
If $\phi $ is a completely positive flow satisfying the quantum stochastic
equation (\ref{2.2}) with $\phi _0\left( B\right) =B$, then the germ-matrix $%
\boldsymbol{\gamma }=\left( \lambda _\nu ^\mu +\imath _\nu ^\mu \right)
_{\nu =+,\bullet }^{\mu =-,\bullet }$ is conditionally completely positive
in the sense 
\begin{equation*}
\sum_{B\in \mathcal{B}}\boldsymbol{\iota }\left( B\right) \boldsymbol{\zeta }%
_B=0\Rightarrow \sum_{B,C\in \mathcal{B}}\langle \boldsymbol{\zeta }_B|%
\boldsymbol{\gamma }\left( B^{*}C\right) \boldsymbol{\zeta }_C\rangle \geq 0.
\end{equation*}
Here $\boldsymbol{\zeta }\in \mathcal{D}\oplus \mathcal{D}_{\bullet },%
\mathcal{D}_{\bullet }=\mathcal{D}\otimes \mathcal{E}$, and $\boldsymbol{%
\iota }=\left( \iota _\nu ^\mu \right) _{\nu =+,\bullet }^{\mu =-,\bullet }$
is the degenerate representation $\iota _\nu ^\mu \left( B\right) =B\delta
_\nu ^{+}\delta _{-}^\mu $, written both with $\boldsymbol{\gamma }$ in the
matrix form as 
\begin{equation}
\boldsymbol{\gamma }=\left( 
\begin{array}{cc}
\gamma & \gamma _{\bullet } \\ 
\gamma ^{\bullet } & \gamma _{\bullet }^{\bullet }%
\end{array}
\right) ,\qquad \text{ }\boldsymbol{\iota }\left( B\right) =\left( 
\begin{array}{cc}
B & 0 \\ 
0 & 0%
\end{array}
\right) ,\qquad  \label{2.5}
\end{equation}
where $\gamma =\lambda _{+}^{-},\quad $ $\gamma ^m=\lambda _{+}^m,\quad $ $%
\gamma _n=\lambda _n^{-},\quad \gamma _n^m=\imath _n^m+\lambda _n^m$ with $%
\imath _n^m\left( B\right) =B\delta _n^m$ such that 
\begin{equation}
\gamma \left( B^{*}\right) =\gamma \left( B\right) ^{*},\qquad \text{ }%
\gamma ^n\left( B^{*}\right) =\gamma _n\left( B\right) ^{*},\qquad \text{ }%
\gamma _n^m\left( B^{*}\right) =\gamma _m^n\left( B\right) ^{*}.  \label{2.6}
\end{equation}
If $\phi $ is subfiltering, then $D=-\lambda _{+}^{-}\left( I\right) $ is a
positive Hermitian form, $\left\langle \eta |D\eta \right\rangle \geq 0$,
for all $\eta \in \mathcal{D}$, and if $\phi $ is contractive, then $%
\boldsymbol{D}=-\boldsymbol{\lambda }\left( I\right) $ is positive in the
sense $\langle \boldsymbol{\eta }|\boldsymbol{D}\boldsymbol{\eta }\rangle
\geq 0$ for all $\boldsymbol{\eta }\in \mathcal{D}\oplus \mathcal{D}%
_{\bullet }$.
\end{theorem}

The proof is given in \cite{Bge, Be97} even for the general (noncommutative)
algebras $\mathfrak{a}$ and $\mathcal{A}$.

Obviously the CCP property for the germ-matrix $\boldsymbol{\gamma }$ is
invariant under the transformation $\boldsymbol{\gamma }\mapsto \boldsymbol{%
\varphi }$ given by 
\begin{equation}
\boldsymbol{\varphi }\left( B\right) =\boldsymbol{\gamma }\left( B\right) +%
\boldsymbol{\iota }\left( B\right) \boldsymbol{K}+\boldsymbol{K}^{*}%
\boldsymbol{\iota }\left( B\right) ,  \label{2.8}
\end{equation}
where $\boldsymbol{K}=\left( K_\nu ^\mu \right) _{\nu =+,\bullet }^{\mu
=-,\bullet }$ is an arbitrary matrix of $K_\nu ^\mu \in \mathcal{L}\left( 
\mathcal{D}\right) $ with $K_{-\nu }^{*\mu }=K_{-\mu }^{\nu *}$. As was
proven in \cite{15, Be96} for the case of finite-dimensional matrix $%
\boldsymbol{\gamma }$ of bounded $\gamma _\nu ^\mu $, see also \cite{LiP},
the matrix elements $K_\nu ^{-}$ can be chosen in such way that the matrix
map $\boldsymbol{\varphi }=\left( \varphi _\nu ^\mu \right) _{\nu =+,\bullet
}^{\mu =-,\bullet }$ becomes CP from $\mathcal{B}$ into the quadratic
matrices of $\varphi _\nu ^\mu \left( B\right) .$ (The other elements can be
chosen arbitrarily, say as $K_{+}^{\bullet }=0$, $K_{\bullet }^{\bullet
}=\frac 12I_{\bullet }^{\bullet }$, because (\ref{2.8}) does not depend on $%
K_{+}^{\bullet },K_{\bullet }^{\bullet }$.) Thus the generator $\boldsymbol{%
\lambda }=\boldsymbol{\gamma }-\boldsymbol{\imath }$ for a quantum
stochastic CP flow $\phi $ can be written (at least in the bounded case) as $%
\boldsymbol{\varphi }-\imath \boldsymbol{K}-\boldsymbol{K}^{*}\imath $: 
\begin{equation}
\lambda _\nu ^\mu \left( B\right) =\varphi _\nu ^\mu \left( B\right)
-B\left( \tfrac 12\delta _\nu ^\mu I+\delta _{-}^\mu K_\nu \right) -\left(
\tfrac 12\delta _\nu ^\mu I+K^\mu \delta _\nu ^{+}\right) B,  \label{2.10}
\end{equation}
where $\varphi _\nu ^\mu :\mathcal{B}\rightarrow \mathcal{B}\left( \mathcal{D%
}\right) $ are matrix elements of the CP map $\boldsymbol{\varphi }$ and $%
K_\nu \in \mathcal{L}\left( \mathcal{D}\right) $, $K^{-}=K_{+}^{*}$, $%
K^m=K_m^{*}$. Now we show that the germ-matrix of this form obeys the CCP
property even in the general case of unbounded $K_\nu ^{-},$ $\varphi _\nu
^\mu \left( B\right) \in \mathcal{B}\left( \mathcal{D}\right) $.

\begin{proposition}
The matrix map $\boldsymbol{\gamma }=\left( \gamma _\nu ^\mu \right) _{\nu
=+,\bullet }^{\mu =-,\bullet }$ given in (\ref{2.8}) by $\quad $%
\begin{equation}
\boldsymbol{\varphi }=\left( 
\begin{array}{cc}
\varphi & \varphi _{\bullet } \\ 
\varphi ^{\bullet } & \varphi _{\bullet }^{\bullet }%
\end{array}
\right) ,\quad \mathrm{and\quad }\boldsymbol{K}=\left( 
\begin{array}{cc}
K & K_{\bullet } \\ 
0 & \frac 12I_{\bullet }^{\bullet }%
\end{array}
\right) ,\;\boldsymbol{K}^{*}=\left( 
\begin{array}{cc}
K^{*} & 0 \\ 
K_{\bullet }^{*} & \frac 12I_{\bullet }^{\bullet }%
\end{array}
\right) ,  \label{2.9}
\end{equation}
with $\varphi =\varphi _{+}^{-},\quad \varphi ^m=\varphi _{+}^m,\quad
\varphi _n=\varphi _n^{-}$ and $\varphi _n^m=\gamma _n^m$ is CCP with
respect to the degenerate representation $\boldsymbol{\iota }=\left( \delta
_{-}^\mu \delta _\nu ^{+}\iota \right) _{\nu =+,\bullet }^{\mu =-,\bullet }$%
, where $\iota \left( B\right) =B$, if $\boldsymbol{\varphi }$ is a CP map.
\end{proposition}

\proof%
If $\boldsymbol{\iota }\left( B_k\right) \boldsymbol{\eta }^k=0$, then 
\begin{equation*}
\langle \boldsymbol{\eta }^k|\boldsymbol{\iota }\left( B_k^{*}B_l\right) 
\boldsymbol{K}+\boldsymbol{K}^{*}\boldsymbol{\iota }\left( B_k^{*}B_l\right) 
\boldsymbol{\eta }^l\rangle
\end{equation*}
\begin{equation*}
=2\func{Re}\left\langle \boldsymbol{\iota }\left( B_k\right) \boldsymbol{%
\eta }^k|\boldsymbol{\iota }\left( B_l\right) \boldsymbol{K}\boldsymbol{\eta 
}^l\right\rangle =0.
\end{equation*}
Hence the CCP for $\boldsymbol{\gamma }$ is equivalent to the CCP property
for (\ref{2.8}) and follows from its CP property: 
\begin{equation*}
\left\langle \boldsymbol{\eta }^k|\boldsymbol{\gamma }\left(
B_k^{*}B_l\right) \boldsymbol{\eta }^l\right\rangle =\left\langle 
\boldsymbol{\eta }^k|\boldsymbol{\varphi }\left( B_k^{*}B_l\right) 
\boldsymbol{\eta }^l\right\rangle \geq 0
\end{equation*}
for such sequences $\boldsymbol{\eta }^k\in \mathcal{D}\oplus \mathcal{D}%
_{\bullet }$.

\section{Construction of quantum CP flows}

The necessary conditions for the stochastic generator $\boldsymbol{\lambda }%
=\left( \lambda _\nu ^\mu \right) _{\nu =+,\bullet }^{\mu =-,\bullet }$ of a
CP flow $\phi $ at $t=0$ are found in the previous section in the form of a
CCP property for the corresponding germ $\boldsymbol{\gamma }=\left( \gamma
_\nu ^\mu \right) _{\nu =+,\bullet }^{\mu =-,\bullet }$. In the next section
we shall show, these conditions are essentially equivalent to the assumption
(\ref{2.10}), corresponding to 
\begin{equation}
\gamma ^m\left( B\right) =\varphi ^m\left( B\right) -K_m^{*}B=\gamma
_m^{*}\left( B\right) ,\quad \gamma \left( B\right) =\varphi \left( B\right)
-K^{*}B-BK,  \label{3.0'}
\end{equation}
where $\boldsymbol{\varphi }=\left( \varphi _\nu ^\mu \right) _{\nu
=+,\bullet }^{\mu =-,\bullet }$ is a CP map with $\varphi _n^m=\gamma _n^m$.
Here we are going to prove under the following conditions for the operators $%
K,K_{\bullet }$ and the maps $\varphi _\nu ^\mu $ that this general form is
also sufficient for the existence of the CP solutions to the quantum
stochastic equation (\ref{2.2}). We are going to construct the minimal
quantum stochastic positive flow $B\mapsto \phi _t\left( B\right) $ for a
given w*-continuous unbounded germ-matrix map of the above form, satisfying
the following conditions.

\begin{enumerate}
\item First, we suppose that the operator $K\in \mathcal{B}\left( \mathcal{D}%
\right) $ generates the one parametric semigroup $\left( e^{-Kt}\right)
_{t>0}$, $e^{-Kr}e^{-Ks}=e^{-K\left( r+s\right) }$ of continuous operators $%
e^{-Kt}\in \mathcal{L}\left( \mathcal{D}\right) $ in the strong sense 
\begin{equation*}
\lim_{t\searrow 0}\frac 1t\left( I-e^{-Kt}\right) \eta =K\eta ,\quad \forall
\eta \in \mathcal{D}.
\end{equation*}
(A contraction semigroup on the Hilbert space $\mathcal{H}$ if $K$ defines
an accretive $K+K^{\dagger }\geq 0$ and so maximal accretive form.)

\item Second, we suppose that the solution $S_t^n,n\in \mathbb{N}$ to the
recurrence 
\begin{equation*}
S_t^{n+1}=S_t^{\circ }-\int_0^tS_{t-r}^{\circ }\sum_{m=1}^\infty K_mS_r^n%
\mathrm{d}\Lambda _{-}^m,\quad S_t^0=S_t^{\circ },
\end{equation*}
where $S_t^{\circ }=e^{-Kt}\otimes T_t\in \mathcal{L}\left( \mathfrak{D}%
\right) $ is the contraction given by the shift co-isometries $T_t:\mathfrak{%
F}\rightarrow \mathfrak{F}$, strongly converges to a continuous operator $%
S_t\in \mathcal{L}\left( \mathfrak{D}\right) $ at $n\longrightarrow \infty $
for each $t>0$.

\item Third, we suppose that the solution $R_t^n,n\in \mathbb{N}$ to the
recurrence 
\begin{equation*}
R_t^{n+1}=S_t^{*}S_t+\int_0^t\mathrm{d}\Lambda _\mu ^\nu \left(
r,S_r^{*}\varphi _\nu ^\mu \left( R_{t-r}^n\right) S_r\right) ,\quad
R_t^0=S_t^{*}S_t,
\end{equation*}
where the quantum stochastic non-adapted integral is understood in the sense 
\cite{20}, weakly converges to a continuous form $R_t\in \mathcal{B}\left( 
\mathfrak{D}\right) $ at $n\longrightarrow \infty $ for each $t>0$.
\end{enumerate}

The first and second assumptions are necessary to define the existence of
free evolution semigroup $S^{\circ }=\left( S_t^{\circ }\right) _{t>0}$ and
its perturbation $S=\left( S_t\right) _{t>0}$ on the product space $%
\mathfrak{D}=\mathcal{D}\otimes \mathfrak{F}$ in the form of multiple
quantum stochastic integral 
\begin{equation}
S_t=S_t^{\circ }+\sum_{n=1}^\infty \left( -1\right) ^n\underset{0<t_1<\ldots
<t_n<t}{\int \cdots \int }K_{m_n}\left( t-t_n\right) \cdot \cdot \cdot
K_{m_1}\left( t_2-t_1\right) S_{t_1}^{\circ }\mathrm{d}\Lambda
_{-}^{m_1}\cdot \cdot \cdot \mathrm{d}\Lambda _{-}^{m_n},  \label{3.5'}
\end{equation}
iterating the quantum stochastic integral equation 
\begin{equation}
S_t=S_t^{\circ }-\int_0^t\sum_{m=1}^\infty K_m\left( t-r\right) S_r\mathrm{d}%
\Lambda _{-}^m,\quad S_0=I,  \label{3.1'}
\end{equation}
where $K_m\left( t\right) =S_t^{\circ }\left( K_m\otimes I\right) $. The
third assumption supplies the weak convergence for the series 
\begin{equation}
R_t=S_t^{*}S_t+\sum_{n=1}^\infty \underset{0<t_1<\ldots <t_n<t}{\int \cdots
\int }\mathrm{d}\Lambda _{\mu _1\ldots \mu _n}^{\nu _1\ldots \nu _n}\left(
t_1,\ldots ,t_n,\varphi _{\nu _1...\nu _n}^{\mu _1...\mu _n}\left(
t_1,\ldots ,t_n,S_{t-t_n}^{*}S_{t-t_n}\right) \right)  \label{3.6'}
\end{equation}
of non-adapted n-tuple integrals, i.e. for the multiple quantum stochastic
integral \cite{20} with 
\begin{equation}
\varphi _{\nu _1\ldots \nu _n}^{\mu _1\ldots \mu _n}\left( t_1,\ldots
,t_n\right) =\varphi _{\nu _1\ldots \nu _{n-1}}^{\mu _1\ldots \mu
_{n-1}}\left( t_1,\ldots ,t_{n-1}\right) \circ \varphi _{\nu _n}^{\mu
_n}\left( t_n-t_{n-1}\right) ,  \label{3.9'}
\end{equation}
where$\quad \varphi _\nu ^\mu \left( t,B\right) =S_t^{*}\varphi _\nu ^\mu
\left( B\right) S_t$.

The following theorem gives a characterization of the evolution semigroup $S$
in terms of cocycles with unbounded coefficients, characterized by Fagnola 
\cite{Fgn} in the isometric and unitary case.

\begin{proposition}
Let the family $V^{\circ }=\left( V_t^{\circ }\right) _{t>0}$ be a quantum
stochastic adapted cocycle, $V_r^{\circ }T_sV_s^{\circ }=T_sV_{r+s}^{\circ }$%
, satisfying the HP differential equation 
\begin{equation}
\mathrm{d}V_t^{\circ }+KV_t^{\circ }\mathrm{d}t+\sum_{m=1}^\infty
K_mV_t^{\circ }\mathrm{d}\Lambda _{-}^m+\sum_{n=1}^\infty V_t^{\circ }%
\mathrm{d}\Lambda _n^n=0,\quad V_0^{\circ }=I.  \label{3.2'}
\end{equation}
Then $S_t=T_tV_t^{\circ }$ is a semigroup solution, $S_rS_s=S_{r+s}$ to the
non-adapted integral equation (\ref{3.1'}) such that $S_t\psi _f=S_t\left(
f^{\bullet }\right) \eta \otimes \delta _{\varnothing },\forall \eta \in 
\mathcal{D}$ on $\psi _f=\eta \otimes f^{\otimes }$ with $f^{\bullet }\in 
\mathfrak{E}_t$. Conversely, if $S=\left( S_t\right) _{t>0}$ is the
non-adapted solution (\ref{3.5'}) to the integral equation (\ref{3.1'}),
then $V_t^{\circ }=T_t^{*}S_t$ is the adapted solution to (\ref{3.2'}),
defined as $V_t^{\circ }\psi _f=S_t\left( f^{\bullet }\right) \eta ,\forall
\eta \in \mathcal{D}$, where $S_t\left( f^{\bullet }\right) =F^{*}S_tF$ is
given by $F\eta =\eta \otimes f^{\otimes }$ with $f^{\bullet }\in \mathfrak{E%
}_t$.
\end{proposition}

\proof%
First let us show that the equation (\ref{3.2'}) is equivalent to the
integral one 
\begin{equation*}
V_t^{\circ }=e^{-Kt}\otimes I_t-\int_0^t\sum_{m=1}^\infty \left( e^{-K\left(
t-r\right) }K_m\otimes I_{t-r}^r\right) V_r^{\circ }\mathrm{d}\Lambda
_{-}^m,\quad V_0^{\circ }=I,
\end{equation*}
where $I_t=T_t^{\dagger }T_t=0^{\Lambda _{\bullet }^{\bullet }\left(
t\right) }$ is the decreasing family of orthoprojectors onto $\mathfrak{F}%
_{[t}$, and $I_s^r=\theta ^r\left( I_s\right) $. Indeed, multiplying both
parts of the integral equation from the left by $\left( e^{K\left(
t-s\right) }\otimes I\right) $ and differentiating the product $e^{K\left(
t-s\right) }V_t^{\circ }$ at $t=s$, we obtain (\ref{3.2'}) by taking into
account that $\mathrm{d}I_t+\sum_{n=1}^\infty I_t\mathrm{d}\Lambda _n^n=0$
and $\mathrm{d}\Lambda _n^n\mathrm{d}\Lambda _{-}^m=0$. Conversely, the
integral equation can be obtained from (\ref{3.2'}) by the integration: 
\begin{eqnarray*}
V_t^{\circ }-e^{-Kt}\otimes I_t &=&\int_0^t\mathrm{d}\left( \left(
e^{-K\left( t-r\right) }\otimes I_{t-r}^r\right) V_r^{\circ }\right) \\
&=&\int_0^t\left( e^{-K\left( t-r\right) }\otimes I_{t-r}^r\right) \left( 
\mathrm{d}V_r^{\circ }+KV_r^{\circ }\mathrm{d}r+V_r^{\circ }\mathrm{d}%
\Lambda _{\bullet }^{\bullet }\right) \\
&=&-\int_0^te^{-K\left( t-r\right) }\left( K_{\bullet }\otimes
I_{t-r}^r\right) V_r^{\circ }\mathrm{d}\Lambda _{-}^{\bullet },
\end{eqnarray*}
where we used that $\mathrm{d}I\left( r\right) =I\left( r\right) \mathrm{d}%
\Lambda _{\bullet }^{\bullet }$ and $\mathrm{d}\left( I\left( r\right)
V_r^{\circ }\right) =\mathrm{d}I\left( r\right) V_r^{\circ }+I\left(
r\right) \mathrm{d}V_r^{\circ }$ for the backward-adapted process $I\left(
r\right) =I_{t-r}^r,\forall r\leq t$. The non-adapted equation (\ref{3.1'})
is obtained by applying the operator $T_t=T_{t-r}T_r$ to both parts of this
integral equation and taking into account the commutativity of $e^{K\left(
r-t\right) }K_m$ with $T_r$. Moreover, due to the adaptiveness of $%
V_t^{\circ }$, $S_t\psi _f=T_t\left( E_tV_t^{\circ }\psi _f\otimes
E_{[t}f^{\otimes }\right) =S_t\left( f^{\bullet }\right) \eta \otimes
f_t^{\otimes }$, where $f_t^{\otimes }=T_tf^{\otimes }$, and $S_t\left(
f^{\bullet }\right) =EV_t^{\circ }F$ is the solution to the equation 
\begin{equation*}
S_t\left( f^{\bullet }\right) =e^{-Kt}+\int_0^te^{-K\left( t-r\right)
}K_{\bullet }f^{\bullet }\left( r\right) S_r\left( f^{\bullet }\right) 
\mathrm{d}r,\ \qquad S_0\left( f^{\bullet }\right) =I.
\end{equation*}
Hence $S_tF=E^{*}S_t\left( f^{\bullet }\right) $ if $f^{\bullet }\in 
\mathfrak{E}_t$, and $F^{*}S_tF=S_t\left( f^{\bullet }\right) $ as $EF=I$.
Since this equation is equivalent to the differential one 
\begin{equation}
\frac{\mathrm{d}}{\mathrm{d}t}S_t\left( f^{\bullet }\right) \eta +\left(
K_{\bullet }f^{\bullet }\left( t\right) +K\right) S_t\left( f^{\bullet
}\right) \eta =0,\quad S_0\left( f^{\bullet }\right) \eta =\eta ,\qquad
\forall \eta \in \mathcal{D},  \label{3.7'}
\end{equation}
the function $t\mapsto S_t\left( f^{\bullet }\right) ,$ $f^{\bullet }\in 
\mathfrak{E}$ is a strongly continuous cocycle, 
\begin{equation*}
S_r\left( f_s^{\bullet }\right) S_s\left( f^{\bullet }\right) =S_{r+s}\left(
f^{\bullet }\right) ,\ \forall r,s>0,\qquad f_s^{\bullet }\left( t\right)
=f^{\bullet }\left( t+s\right) ,\quad S_0\left( f^{\bullet }\right) =I.
\end{equation*}
As was proved in \cite{20}, the multiple integral (\ref{3.5'}) gives a
solution to the integral equation (\ref{3.1'}), and so the multiple integral
for $V_t^{\circ }\psi _f=S_t\left( f^{\bullet }\right) \eta \otimes
f^{\otimes },$ 
\begin{equation*}
S_t\left( f^{\bullet }\right) =e^{-Kt}+\sum_{n=1}^\infty \left( -1\right) ^n%
\underset{0<t_1<\ldots <t_n<t}{\int \cdots \int }K\left( t,t_n\right) \cdot
\cdot \cdot K\left( t_2,t_1\right) e^{-Kt_1}\mathrm{d}t_1\cdot \cdot \cdot 
\mathrm{d}t_n,
\end{equation*}
where $K\left( t,r\right) =e^{-K\left( t-r\right) }K_{\bullet }f^{\bullet
}\left( r\right) $, corresponding to the iteration of the integral equation
for $V_t^{\circ }$ on $\psi _f$, satisfies the HP equation (\ref{3.2'}). 
\endproof%

The following theorem reduces the problem of solving of differential
evolution equations to the problem of iteration of integral equations
similar to the nonstochastic case \cite{Be89, Chb}.

\begin{proposition}
Let $S_t=T_tV_t^{\circ }$, where $V_t^{\circ }\in \mathcal{L}\left( 
\mathfrak{D}\right) $ are continuous operators defining the adapted cocycle
solution to the equation (\ref{3.2'}). Then the linear stochastic evolution
equation (\ref{2.2}) is equivalent to the quantum non-adapted (in the sense
of \cite{20}) integral equation 
\begin{equation}
\phi _t\left( B\right) =S_t^{*}BS_t+\int_0^t\mathrm{d}\Lambda _\mu ^\nu
\left( r,\phi _r\left[ \varphi _\nu ^\mu \left( S_{t-r}^{*}BS_{t-r}\right) %
\right] \right)  \label{3.3'}
\end{equation}
with $\phi _0\left( B\right) =B\in \mathcal{B}$, where $\varphi _\nu ^\mu $
are extended onto $\mathfrak{B}$ in the normal way by w*-continuity and
linearity as $\varphi _\nu ^\mu \left( B\otimes Z\right) =\overline{\varphi }%
_\nu ^\mu \left( B\right) \otimes Z$ for $B\in \overline{\mathcal{B}}$, $%
Z\in \mathcal{B}\left( \mathfrak{F}\right) $.
\end{proposition}

The proof is given in \cite{Bpfg, Be97}.

\begin{theorem}
Let $\boldsymbol{\varphi }$ be a w*-continuous CP-map, and $%
S_t=T_tV_t^{\circ } $ be given by the solution to the quantum stochastic
equation (\ref{3.2'}). Then the solutions to the evolution equation (\ref%
{2.2}) with the generators, corresponding to (\ref{3.0'}), have the CP
property, and satisfy the submartingale (contractivity) condition $\phi
_t\left( I\right) \leq \epsilon _t\left[ \phi _s\left( I\right) \right] $
for all $t<s$ if $\varphi \left( I\right) \leq K+K^{\dagger }$ ($\phi
_t(I)\leq \phi _s(I)$ if $\boldsymbol{\varphi }(I)\leq \boldsymbol{K}+%
\boldsymbol{K}^{\dagger }$). The minimal solution can be constructed in the
form of multiple quantum stochastic integral in the sense \cite{20} as the
series 
\begin{equation}
\phi _t\left( B\right) =\sum_{n=0}^\infty \underset{0<t_1<\ldots <t_n<t}{%
\int \cdots \int }\mathrm{d}\Lambda _{\mu _1\ldots \mu _n}^{\nu _1\ldots \nu
_n}\left( t_1,\ldots ,t_n,\varphi _{\nu _1\ldots \nu _n}^{\mu _1\ldots \mu
_n}\left( t_1,\ldots ,t_n,S_{t-t_n}^{*}BS_{t-t_n}\right) \right)
\label{3.8'}
\end{equation}
of non-adapted n-tuple CP integrals with $S_t^{*}BS_t$ at $n=0$ and 
\begin{equation*}
\varphi _{\nu _1\ldots \nu _n}^{\mu _1\ldots \mu _n}\left( t_1,\ldots
,t_n\right) =\varphi _{\nu _1}^{\mu _1}\left( t_1\right) \circ \varphi _{\nu
_2}^{\mu _2}\left( t_2-t_1\right) \circ \ldots \circ \varphi _{\nu _n}^{\mu
_n}\left( t_n-t_{n-1}\right) ,
\end{equation*}
where $\varphi _\nu ^\mu \left( t,B\right) =S_t^{*}\varphi _\nu ^\mu \left(
B\right) S_t$. If $\boldsymbol{\varphi }$ is bounded, then the solution to
the equation is unique, and $\phi _t\left( I\right) =\epsilon _t\left[ \phi
_s\left( I\right) \right] $ for all $t<s$ if $K+K^{\dagger }=\varphi \left(
I\right) $ ($\phi _t(I)=I$ if $\boldsymbol{K}+\boldsymbol{K}^{\dagger }=%
\boldsymbol{\varphi }(I)$).
\end{theorem}

The proof is given in \cite{Bpfg, Be97}.

\section{The structure of the generators and flows}

First, let us prove the structure (\ref{2.10}) for the (unbounded)
form-generator of CP flows over the algebra $\mathcal{B}=\mathcal{L}\left( 
\mathcal{H}\right) $ of all bounded operators, assuming that $\mathcal{A}=0$%
. This algebra contains the one-dimensional operators $|\eta ^{\prime
}\rangle \langle \eta ^0|:\eta \mapsto \left\langle \eta ^0|\eta
\right\rangle \eta ^{\prime }$ given by the vectors $\eta ^0,\eta ^{\prime
}\in \mathcal{H}$.

Let us fix a vector $\boldsymbol{\eta }^0\in \mathcal{D}\oplus \mathcal{D}%
_{\bullet }$ with the unit projection $\eta ^0\in \mathcal{D}$, $\left\|
\eta ^0\right\| =1$, and make the following assumption of the weak
continuity for the linear operator $\eta ^{\prime }\mapsto $ $\boldsymbol{%
\gamma }\left( |\eta ^{\prime }\rangle \langle \eta ^0|\right) \boldsymbol{%
\eta }^0$.

\begin{enumerate}
\item[0)] The sequence $\boldsymbol{\eta }_n^{\prime }=\boldsymbol{\gamma }%
\left( |\eta _n^{\prime }\rangle \langle \eta ^0|\right) \boldsymbol{\eta }%
^0\in \mathcal{D}^{\prime }\oplus \mathcal{D}_{\bullet }^{\prime }$ of
anti-linear forms 
\begin{equation*}
\boldsymbol{\eta }\in \mathcal{D}\oplus \mathcal{D}_{\bullet }\mapsto
\left\langle \boldsymbol{\eta }|\boldsymbol{\eta }_n^{\prime }\right\rangle
:=\left\langle \boldsymbol{\eta }|\boldsymbol{\gamma }\left( |\eta
_n^{\prime }\rangle \langle \eta ^0|\right) \boldsymbol{\eta }^0\right\rangle
\end{equation*}
converges for each sequence $\eta _n^{\prime }\in \mathcal{H}$ converging in 
$\mathcal{D}^{\prime }\supseteq \mathcal{H}$.
\end{enumerate}

\begin{proposition}
Let the CCP germ-matrix $\boldsymbol{\gamma }$ satisfy the above continuity
condition for a given $\boldsymbol{\eta }^0$. Then there exist strongly
continuous operators $K\in \mathcal{L}\left( \mathcal{D}\right) ,K_{\bullet
}:\mathcal{D}_{\bullet }\rightarrow \mathcal{D}$ defining the matrix
operator $\boldsymbol{K}$ in (\ref{2.10}), such that the matrix map (\ref%
{2.8}) is CP, and there exists a Hilbert space $\mathcal{K}$, a $*$%
-representation $\jmath :B\mapsto B\otimes J$ of $\mathcal{B}=\mathcal{L}%
\left( \mathcal{H}\right) $ on the Hilbert product $\mathcal{G}=\mathcal{H}%
\otimes \mathcal{K}$, given by an orthoprojector $J$ in $\mathcal{K}$, such
that 
\begin{equation}
\boldsymbol{\varphi }\left( B\right) =\left( L^\mu \jmath \left( B\right)
L_\nu \right) _{\nu =+,\bullet }^{\mu =-,\bullet }=\boldsymbol{L}^{*}\jmath
\left( B\right) \boldsymbol{L}.  \label{4.0}
\end{equation}
Here $\boldsymbol{L}=\left( L,L_{\bullet }\right) $ is a strongly continuous
operator $\mathcal{D}\oplus \mathcal{D}_{\bullet }\rightarrow \mathcal{G}$
with $L=L_{+}$, $L^{-}=L^{*}$, $L^{\bullet }=L_{\bullet }^{*}$ which is
always possible to make 
\begin{equation}
\left\langle \eta ^0\otimes e|\boldsymbol{L}\boldsymbol{\eta }%
^0\right\rangle =0,\quad \quad \forall e\in \mathcal{K}_1,  \label{4.1}
\end{equation}
where $\mathcal{K}_1=J\mathcal{K}$. If $D=-\lambda \left( I\right) \geq 0$,
then one can make $L^{*}L=K+K^{\dagger }$ in a canonical way $L=L^{\circ }$,
and in addition one can make $L^{*}L_{\bullet }=K_{\bullet }$, $L_{\bullet
}^{*}L_{\bullet }=I_{\bullet }^{\bullet }$, where $I_{\bullet }^{\bullet
}=I\delta _{\bullet }^{\bullet }$ for a canonical $L_{\bullet }=L_{\bullet
}^{\circ }$ if $\boldsymbol{D}=-\boldsymbol{\lambda }\left( I\right) \geq 0$.
\end{proposition}

The proof is given in \cite{Be97}.

Thus we have proved that the equation (\ref{2.2}) for a completely positive
quantum stochastic flows over $\mathcal{B}=\mathcal{L}\left( \mathcal{H}%
\right) $ has the following general form 
\begin{equation*}
\mathrm{d}\phi _{t}\left( B\right) +\phi _{t}\left( K^{\ast }B+BK-L^{\ast
}\jmath \left( B\right) L\right) \mathrm{d}t=\sum_{m,n=1}^{\infty }\phi
_{t}\left( L_{m}^{\ast }\jmath \left( B\right) L_{n}-B\delta _{n}^{m}\right) 
\mathrm{d}\Lambda _{m}^{n}
\end{equation*}
\begin{equation*}
+\sum_{m=1}^{\infty }\phi _{t}\left( L_{m}^{\ast }\jmath \left( B\right)
L-K_{m}^{\ast }B\right) \mathrm{d}\Lambda _{m}^{+}+\sum_{n=1}^{\infty }\phi
_{t}\left( L^{\ast }\jmath \left( B\right) L_{n}-BK_{n}\right) \mathrm{d}%
\Lambda _{-}^{n},
\end{equation*}
generalizing the Lindblad form \cite{19} for the semigroups of completely
positive maps. This can be written in the tensor notation form as 
\begin{equation}
\mathrm{d}\phi _{t}\left( B\right) =\phi _{t}\left( L_{\alpha }^{\star \mu
}\jmath _{\beta }^{\alpha }\left( B\right) L_{\nu }^{\beta }-\imath _{\nu
}^{\mu }\left( B\right) \right) \mathrm{d}\Lambda _{\mu }^{\nu }=\phi
_{t}\left( \mathbf{L}^{\star }\text{\textbf{\j }}\left( B\right) \mathbf{L}-%
\text{\textbf{\i }}\left( B\right) \right) \cdot \mathrm{d}\boldsymbol{%
\Lambda },  \label{4.3}
\end{equation}
where the summation is taken over all $\alpha ,\beta =-,\circ ,+$ and $\mu
,\nu =-,\bullet ,+$, $\jmath _{-}^{-}\left( B\right) =B=\jmath
_{+}^{+}\left( B\right) $, $\jmath _{\circ }^{\circ }\left( B\right) =\jmath
\left( B\right) $, $\jmath _{\beta }^{\alpha }\left( B\right) =0$ if $\alpha
\neq \beta $, $\imath _{\nu }^{\mu }\left( B\right) =B\delta _{\nu }^{\mu }$%
, and $\mathbf{L}^{\star }\mathbf{=}\left[ L_{\beta }^{\star \mu }\right]
_{\beta =-,\circ ,+}^{\mu =-,\bullet ,+}$ with $L_{-\alpha }^{\star \mu
}=L_{-\mu }^{\alpha \ast }$ is the triangular matrix, pseudoadjoint to $%
\mathbf{L=}\left[ L_{\nu }^{\alpha }\right] _{\nu =-,\bullet ,+}^{\alpha
=-,\circ ,+}$ with $L_{-}^{-}=I=L_{+}^{+}$, 
\begin{equation*}
L_{\bullet }^{\circ }=L_{\bullet },\qquad L_{+}^{\circ }=L,\qquad L_{\bullet
}^{-}=-K_{\bullet },\qquad L_{+}^{-}=-K.
\end{equation*}
(All other $L_{\nu }^{\alpha }$ are zero.) If the Hilbert space $\mathcal{K}$
is separable, $\mathcal{K}_{1}=\ell ^{2}\left( \mathbb{N}_{1}\right) $ for a
subset $\mathbb{N}_{1}\subseteq \mathbb{N}$. Then the equation (\ref{4.3})
can be resolved as $\phi _{t}\left( B\right) =V_{t}^{\ast }\left( B\otimes
I_{t}\right) V_{t}$, where $V=\left( V_{t}\right) _{t>0}$ is an (unbounded)
cocycle on the product $\mathcal{D}\otimes \mathfrak{F}$ with Fock space $%
\mathfrak{F}$ over the Hilbert space $L^{2}\left( \mathbb{N\times R}%
_{+}\right) $ of the quantum noise, and $I_{t}$ is a decreasing family of
orthoprojectors in $\mathfrak{F}$, satisfying the stochastic equation $%
\mathrm{d}I_{t}+\sum_{n\notin \mathbb{N}_{1}}I_{t}\mathrm{d}\Lambda
_{n}^{n}=0$ with $I_{0}=I$. The cocycle $V$ can be found from the quantum
stochastic equation $\mathrm{d}V_{t}=(L_{\nu }^{\mu }-I\delta _{\nu }^{\mu
})V_{t}\mathrm{d}\Lambda _{\mu }^{\nu }$ with $V_{0}=I\otimes I$ of the form 
\begin{equation}
\mathrm{d}V_{t}+KV_{t}\mathrm{d}t+\sum_{n=1}^{\infty }K_{n}V_{t}\mathrm{d}%
\Lambda _{-}^{n}=\sum_{m,n=1}^{\infty }\left( L_{n}^{m}-I\delta
_{n}^{m}\right) V_{t}\mathrm{d}\Lambda _{m}^{n}+\sum_{m=1}^{\infty
}L^{m}V_{t}\mathrm{d}\Lambda _{m}^{+},\qquad  \label{4.4}
\end{equation}
where $L_{n}^{m\text{ }}$ and $L^{m}$ are the operators in $\mathcal{D}$,
defining 
\begin{eqnarray}
\varphi _{n}^{m}\left( B\right) &=&\sum_{k\in \mathbb{N}_{1}}L_{m}^{k\ast
}BL_{n}^{k},\qquad \varphi \left( B\right) =\sum_{k\in \mathbb{N}%
_{1}}L^{k\ast }BL^{k}  \label{4.5} \\
\varphi ^{m}\left( B\right) &=&\sum_{k\in \mathbb{N}_{1}}L_{m}^{k\ast
}BL^{k},\qquad \varphi _{n}\left( B\right) =\sum_{k\in \mathbb{N}%
_{1}}L^{k\ast }BL_{n}^{k}\qquad  \notag
\end{eqnarray}
with $\sum_{k=1}^{\infty }L^{k\ast }L^{k}=K+K^{\dagger }$ if $D\geq 0$, and
in addition $\sum_{k=1}^{\infty }L^{k\ast }L_{n}^{k}=K_{n}$, $%
\sum_{k=1}^{\infty }L_{m}^{k\ast }L_{n}^{k}=I\delta _{n}^{m}$ if $%
\boldsymbol{D}\geq 0$. \allowbreak The formal derivation of the equation (%
\ref{4.4}) from (\ref{4.3}) is obtained by a simple application of the HP It%
\^{o} formula. The martingale $M_{t}$, describing the density operator for
the output state of $\Lambda \left( t,a\right) $, is then defined as $%
M_{t}=V_{t}^{\ast }V_{t} $.

The following theorem ensures the existence of a $*$-representation $\iota
:\Lambda \left( t,a\right) \mapsto \Lambda \left( t,i\left( a\right) \right)
:=i_\beta ^\alpha \left( a\right) \Lambda _\alpha ^\beta \left( t\right) $
of the quantum stochastic process (\ref{0.3}), commuting with $Y_t=\phi
_t\left( B\right) $ for all $a\in \mathfrak{a},B\in \mathcal{L}\left( 
\mathcal{H}\right) $, with $A=X\left( 0\right) =0$, in the form 
\begin{equation*}
\Lambda \left( t,i\left( a\right) \right) =i_{\circ }^{\circ }\left(
a\right) \Lambda _{\circ }^{\circ }\left( t\right) +i_{+}^{\circ }\left(
a\right) \Lambda _{\circ }^{+}\left( t\right) +i_{\circ }^{-}\left( a\right)
\Lambda _{-}^{\circ }\left( t\right) +i_{+}^{-}\left( a\right) \Lambda
_{-}^{+}\left( t\right) \text{.}
\end{equation*}
Here $\boldsymbol{i}\mathbf{=}\left( i_\beta ^\alpha \right) _{\beta
=+,\circ }^{\alpha =-,\circ }$ is a $\star $-representation 
\begin{equation*}
i_\beta ^\alpha \left( a^{\star }a\right) =i_{\circ }^\alpha \left( a^{\star
}\right) i_\beta ^{\circ }\left( a\right) ,\quad i_{-\beta }^\alpha \left(
a^{\star }\right) =i_{-\alpha }^\beta \left( a\right) ^{*}
\end{equation*}
of the It\^{o} algebra $\mathfrak{a}$ in the operators $i_\beta ^\alpha
\left( a\right) :\mathcal{K}_\beta \rightarrow \mathcal{K}_\alpha $, with a
domain $\mathcal{K}_{\circ }\subseteq \mathcal{K}$, $\mathcal{K}_{-}=\mathbb{%
C=}\mathcal{K}_{+}$, and $\Lambda _\alpha ^\beta \left( t\right) $ are the
canonical quantum stochastic integrators in the Fock space $\Gamma \left( 
\mathfrak{K}\right) $ over $\mathfrak{K}=L_{\mathcal{K}}^2\left( \mathbb{R}%
_{+}\right) $, the space of $\mathcal{K}$-valued square-integrable functions
on $\mathbb{R}_{+}$.We shall extend $\boldsymbol{i}$ to the triangular
matrix representation $\mathbf{i=}\left[ i_\beta ^\alpha \right] _{\beta
=-,\circ ,+}^{\alpha =-,\circ ,+}$ on the pseudo-Hilbert space $\mathbb{%
C\oplus }\mathcal{K}\oplus \mathbb{C}$ with the Minkowski metrics tensor $%
\mathbf{g=}\left[ \delta _{-\beta }^\alpha \right] =\mathbf{g}^{-1}$, by $%
i_\beta ^{+}\left( a\right) =0=i_{-}^\alpha \left( a\right) $, for all $a\in 
\mathfrak{a} $, as it was done for $\mathbf{a=}\left[ a_\nu ^\mu \right]
_{\nu =-,\bullet ,+}^{\mu =-,\bullet ,+}$, and denote the ampliation $%
I\otimes i_\beta ^\alpha \left( a\right) $ again as $i_\beta ^\alpha \left(
a\right) $. Note that if the stochastic generator of the form (\ref{2.10})
is restricted onto an operator algebra $\mathcal{B}\subseteq \mathcal{L}%
\left( \mathcal{H}\right) $ with the weak closure $\mathcal{\bar{B}}=%
\mathcal{A}^c$, and all the sesquilinear forms $\gamma _\nu ^\mu \left(
B\right) $, $B\in \mathcal{B} $ commute with the $*$-algebra $\mathcal{A}%
\subset \mathcal{L}\left( \mathcal{D}\right) $, then $\lambda _\nu ^\mu
\left( B\right) \in \mathcal{\bar{B}}$.

\begin{proposition}
Let $\boldsymbol{b}=\boldsymbol{\gamma }\left( B\right) -\boldsymbol{\imath }%
\left( B\right) $ satisfy the commutativity conditions (\ref{2.4}) for all $%
a\in \mathfrak{a}$, $B\in \mathcal{L}\left( \mathcal{H}\right) $. Then there
exists a $\star $-representation $a\mapsto \boldsymbol{i}\left( a\right) $
of the It\^{o} algebra $\mathfrak{a}$, defining the operators $i_\beta
^\alpha \left( a\right) :\mathcal{K}_\beta \rightarrow \mathcal{K}_\alpha $,
with $i_\beta ^\alpha \left( a\right) ^{*}\mathcal{K}_\alpha \subseteq 
\mathcal{K}_\beta $, where $\mathcal{K}_{-}=\mathbb{C=}\mathcal{K}_{+}$,
such that $L_\mu ^\alpha \left( I\otimes a_\nu ^\mu \right) =\left( I\otimes
i_\beta ^\alpha \left( a\right) \right) L_\nu ^\beta $ for all $a\in 
\mathfrak{a}.$ By omitting $I\otimes $ this can be written as 
\begin{eqnarray}
L_{\bullet }a_{\bullet }^{\bullet } &=&i\left( a\right) L_{\bullet },\quad
\quad a_{+}^{-}-K_{\bullet }a_{+}^{\bullet }=i^{-}\left( a\right)
L+i_{+}^{-}\left( a\right) ,  \label{4.6} \\
L_{\bullet }a_{+}^{\bullet } &=&i\left( a\right) L+i_{+}\left( a\right)
,\quad \quad a_{\bullet }^{-}-K_{\bullet }a_{\bullet }^{\bullet
}=i^{-}\left( a\right) L_{\bullet },  \notag
\end{eqnarray}
where we take the convention $i^{-}=i_{\circ }^{-}$, $i_{+}=i_{+}^{\circ }$
and $i=i_{\circ }^{\circ }$. If $\left[ A,\gamma _\nu ^\mu \left( B\right) %
\right] =0$ for all $A\in \mathcal{A}$ and $B\in \mathcal{B}$, where $%
\mathcal{B}\subseteq \mathcal{L}\left( \mathcal{H}\right) $ is a $*$%
-subalgebra, and $\mathcal{\bar{B}}=\mathcal{A}^c$, then there exists a
triangular $\star $-representation $\mathbf{j=}\left[ j_\beta ^\alpha \right]
_{\beta =-,\circ ,+}^{\alpha =-,\circ ,+}$ of the operator algebra $\mathcal{%
A}$ with $j_{\circ }^{\circ }\left( I\right) =J$ such that 
\begin{equation}
\mathbf{J\QTR{mathbf}{LA}}=\mathbf{j}\left( A\right) \mathbf{\QTR{mathbf}{%
L,\quad }}\left[ \mathbf{j}\left( A\right) ,\mathbf{i}\left( a\right) \right]
=0\mathbf{,\quad }\left[ \mathbf{j}\left( A\right) ,\text{\emph{\textbf{\j }}%
}\left( B\right) \right] =0,\quad \forall A\in \mathcal{A},a\in \mathfrak{a}%
,B\in \mathcal{B}.  \label{4.7}
\end{equation}
\end{proposition}

The proof is given in \cite{Be97} even for the general (noncommutative) It%
\^{o} algebra $\mathfrak{a}$.

Now we are going to construct the quantum stochastic dilation for the flow $%
\phi _t\left( B\right) $ and the quantum state generating function $%
\vartheta _t^a=\epsilon \left[ R_tW\left( t,a\right) \right] $ of the output
process $\boldsymbol{\Lambda }\left( t,a\right) $ in the form 
\begin{equation*}
\phi _t\left( B\right) =V_t^{*}\left( I_t\otimes B\right) V_t,\quad
\vartheta _t\left( g\right) =\epsilon \left[ V_t^{*}\left( W_t^a\otimes
I\right) V_t\right] ,\quad \forall B\in \mathcal{L}\left( \mathcal{H}\right)
,a\in \mathfrak{a},
\end{equation*}
where $V_t$ is an operator on $\mathfrak{D}$ into $\Gamma \left( \mathfrak{K}%
\right) \otimes \mathfrak{D}$, intertwining the Weyl operators $W\left(
t,a\right) $ with the operators $W_t^a=W\left( t,i\left( a\right) \right)
I_t $ in the Fock space $\Gamma \left( \mathfrak{K}\right) $, 
\begin{equation*}
\mathrm{d}W\left( t,i\left( a\right) \right) =W\left( t,i\left( a\right)
\right) \mathrm{d}\Lambda \left( t,i\left( a\right) \right) ,\quad W\left(
0,i\left( a\right) \right) =I,
\end{equation*}
and $I_t\geq I_s,\forall t\leq s$ is a decreasing family of orthoprojectors.

In order to prove the existence of the Fock space dilation, we need the
following assumptions in addition to the continuity assumptions of this and
previous sections.

\begin{enumerate}
\item[1)] The minimal quantum stochastic $*$-flow \cite{EvH} over the
operator algebra $\mathcal{A}$, resolving the quantum Langevin equation 
\begin{equation*}
\mathrm{d}\tau _t\left( A\right) =\tau _t\left( \mathbf{j}\left( A\right) -%
\text{\textbf{\i }}\left( A\right) \right) \cdot \mathrm{d}\boldsymbol{%
\Lambda },\quad A\in \mathcal{A},
\end{equation*}
where $\mathbf{j}\left( I\right) =\mathbf{J}\otimes I$, \textbf{\i }$(A)=%
\mathbf{I}\otimes A$ , constructed by its iterations with $\tau _0\left(
A\right) =I\otimes A$ as it was done in the Sec 4 for the flow $\phi $, is
the multiplicative flow, satisfying the condition $\tau _t\left( I\right)
=I_t\otimes I$, where $I_t$ is the solution to $\mathrm{d}I_t=\left(
J-I\right) _{\circ }^{\circ }I_t\mathrm{d}\Lambda _{\circ }^{\circ }$ with $%
I_0=I$.

\item[2)] The operators $L_\nu \left( \bar{e}\right) =\left( I\otimes
e^{*}\right) L_\nu $, given for all $e\in \mathcal{G}$ as $\left\langle
L_\nu \left( \bar{e}\right) \eta |\eta ^{\prime }\right\rangle =\left\langle
L_\nu \eta |\eta ^{\prime }\otimes e\right\rangle \quad \forall \eta \in 
\mathcal{D},\eta ^{\prime }\in \mathcal{D}^{\prime }$, are strongly
continuous onto $\mathcal{D}$. This is necessary for the weak definition of
the operators $V_t\left( \sigma \right) :\mathfrak{D}\rightarrow \mathcal{K}%
^{\otimes \left| \sigma \right| }\otimes \mathfrak{D}_{[t}$ on finite
subsets $\sigma \subset [0,t)$ by the recurrence 
\begin{equation*}
V_t\left( \sigma \right) \psi =\left( I^{\otimes \left| \sigma \right|
}\otimes V_t^{\circ }\left( s\right) \right) \left( LV_s\left( \sigma
\backslash s\right) \psi +\sum_{n=1}^\infty L_nV_s\left( \sigma \backslash
s\right) \psi ^n\left( s\right) \right) ,
\end{equation*}
with $V_t(\emptyset )=V_t^{\circ }\left( 0\right) ,\quad s=\max \sigma $.
Here $V_t^{\circ }\left( s\right) =T_s^{*}V_{t-s}^{\circ }T_s$, $V_t^{\circ
} $ is the solution to the equation (\ref{3.1'}), the operators $L_\nu :%
\mathcal{D}\rightarrow \mathcal{K}\otimes \mathcal{D}$ act on $\mathcal{K}%
^{\otimes \left| \sigma \backslash s\right| }\otimes \mathfrak{D}_{[s}$ as $%
I^{\otimes \left| \sigma \backslash s\right| }\otimes L_\nu \otimes I_{[s}$,
and $\psi ^n\left( s\right) $ are the components of $\psi ^{\bullet }\left(
\tau ,s\right) =\psi \left( \tau \sqcup s\right) $, where $\tau \sqcup s$ is
defined for almost all $s$ ($s\notin \sigma $) as the disjoint union of the
single point $\left\{ s\right\} $ with a finite subset $\tau \in \mathbb{R}%
_{+}$.

\item[3)] The operator-valued function $\sigma \mapsto $ $V_t\left( \sigma
\right) $, defined for all such $\sigma \in \Gamma _t$, is weakly square
integrable for each $t$ with respect to the measure $\mathrm{d}\sigma
=\prod_{s\in \sigma }\mathrm{d}s$ in the sense 
\begin{equation*}
\int_{\Gamma _t}\left\| V_t\left( \sigma \right) \psi \right\| ^2\mathrm{d}%
\sigma :=\sum_{n=0}^\infty \underset{0<s_1\ldots s_n<t}{\int \ldots \int }%
\left\| V_t\left( s_{1,\ldots ,}s_n\right) \psi \right\| ^2\mathrm{d}%
s_1\ldots \mathrm{d}s_1<\infty ,
\end{equation*}
for all $\psi \in \mathfrak{D}.$ Thus the operators $V_t\left( \cdot \right) 
$ define a Fock space one $V_t:$ $\mathfrak{D}\rightarrow \Gamma \left( 
\mathfrak{K}_t\right) \otimes \mathfrak{D}_{[t}$. They form a cocycle, $%
V_{t-r}^r\left( \sigma \right) V_r\left( \sigma \right) =V_t\left( \sigma
\right) $, where $V_s^r\left( \sigma \right) =I^{\otimes \left| \sigma
_r\right| }\otimes T_r^{*}V_s\left( \sigma _{[r}-r\right) T_r$.
\end{enumerate}

\begin{theorem}
Under the given assumptions 0)--3) there exist:

\begin{enumerate}
\item[(i)] A cocycle dilation $V_t:\mathfrak{D}\rightarrow \Gamma \left( 
\mathfrak{K}_t\right) \otimes \mathfrak{D}$ of the minimal CP flow $\phi $,
intertwining the Weyl operator $W\left( t,a\right) $ with $W_t^a$: 
\begin{equation}
V_t\left( I\otimes W\left( t,a\right) \right) =\left( W_t^a\otimes I\right)
V_t,\quad \phi _t\left( B\right) =V_t^{*}\left( I_t\otimes B\right) V_t\text{
},\quad \forall a\in \mathfrak{a},B\in \mathcal{L}\left( \mathcal{H}\right) ,
\label{4.8}
\end{equation}
where $I_t\leq I_s$, $\forall t<s$ are orthoprojectors in $\Gamma \left( 
\mathfrak{K}\right) $;

\item[(ii)] A $*$-multiplicative flow $\tau =\left( \tau _t\right) $ over $%
\mathcal{A}$ in $\Gamma \left( \mathfrak{K}\right) \otimes \mathcal{H}$ with
the properties $\tau _t\left( I\right) =I_t$, 
\begin{equation}
V_tA=\tau _t\left( A\right) V_t,\quad \left[ \tau _t\left( A\right) ,W_t^a%
\right] =0,\quad \left[ \tau _t\left( A\right) ,I\otimes B\right] =0,\quad
\forall A\in \mathcal{A},a\in \mathfrak{a},B\in \mathcal{B}\text{.}
\label{4.9}
\end{equation}

\item[(iii)] If $\lambda \left( I\right) \leq 0$, then one can make $%
M_t=V_t^{*}V_t$ martingale, and, if $\boldsymbol{\lambda }\left( I\right)
\leq 0$, one can make $V_t$ isometric, $V_t^{*}V_t=I$.

\item[(iv)] Moreover, let $U=\left( U_t\right) _{t\geq 0}$ be a
one-parametric weakly continuous cocycle of unitary operators on $\Gamma
\left( \mathfrak{K}\right) \otimes \mathcal{H}\otimes \Gamma \left( 
\mathfrak{E}\right) $ , satisfying the quantum stochastic equation 
\begin{eqnarray}
&&\mathrm{d}U_t+\left( K\mathrm{d}t+K_{\bullet }^{-}\mathrm{d}\Lambda
_{-}^{\bullet }+K_{\circ }^{-}\mathrm{d}\Lambda _{-}^{\circ }\right) U_t
\label{4.10} \\
&=&\left( L_{+}^{\circ }\mathrm{d}\Lambda _{\circ }^{+}-I_{\bullet
}^{\bullet }\mathrm{d}\Lambda _{\bullet }^{\bullet }+J_{\bullet }^{\circ }%
\mathrm{d}\Lambda _{\circ }^{\bullet }+J_{\circ }^{\bullet }\mathrm{d}%
\Lambda _{\bullet }^{\circ }+\left( J_{\circ }^{\circ }-I_{\circ }^{\circ
}\right) \mathrm{d}\Lambda _{\circ }^{\circ }\right) U_t  \notag
\end{eqnarray}
with $U_0=I$ and the necessary differential unitarity conditions 
\begin{equation*}
K+K^{\dagger }=L_{\circ }^{-}L_{+}^{\circ },\ K_{\bullet }^{-}=L_{\circ
}^{-}J_{\bullet }^{\circ },\ J_{\circ }^{\bullet }J_{\bullet }^{\circ
}=I_{\bullet }^{\bullet },\ K_{\circ }^{-}=L_{\circ }^{-}J_{\circ }^{\circ
},\ J_{\circ }^{\circ }=I_{\circ }^{\circ }-J_{\bullet }^{\circ }J_{\circ
}^{\bullet },
\end{equation*}
where $L_{\circ }^{-}=L_{+}^{\circ *}$, $J_{\circ }^{\bullet }=J_{\bullet
}^{\circ *}$. If $\lambda \left( I\right) \leq 0$ and $L_{+}^{\circ
}=L^{\circ }$ is the canonical operator in the dilation (\ref{4.0}), then 
\begin{equation}
\left\langle \psi |\left( A\otimes I\right) \phi _t^a\left( B\right) \psi
\right\rangle =\left\langle U_t\left( \delta _\emptyset \otimes \psi \right)
|\left( \tau _t^a\left( A\right) \left( I\otimes B\right) \right) U_t\left(
\delta _\emptyset \otimes \psi \right) \right\rangle  \label{4.11}
\end{equation}
for all $A\in \mathcal{A},a\in \mathfrak{a},B\in \mathcal{B}$ and any
initial $\psi =\eta \otimes \delta _\emptyset $, $\eta \in \mathcal{D}$,
where 
\begin{equation*}
\phi _t^a\left( B\right) =\left( I\otimes W\left( t,a\right) \right) \phi
_t\left( B\right) ,\quad \tau _t^a\left( A\right) =\left( W_t^a\otimes
I\right) \tau _t\left( A\right) .
\end{equation*}
This unitary cocycle dilation is valied for any state $\psi \in \mathcal{D}%
\otimes \mathfrak{F}$ if in addition $J_{\bullet }^{\circ }=L_{\bullet
}^{\circ } $ is the canonical isometry in (\ref{4.0}) for the case $%
\boldsymbol{\lambda }\left( I\right) \leq 0$.
\end{enumerate}
\end{theorem}

\proof%
(Sketch). The cocycle $V=\left( V_t\right) _{t>0}$ is recurrently
constructed due to the above assumptions (1)--(3). It obviously intertwines
the Weyl operators (\ref{2.1b}) with the operators $W_t^a$, acting in the
same way in $\Gamma \left( \mathfrak{K}\right) $, by virtue of the property (%
\ref{4.6})$.$

Let us denote by $\mathfrak{K}_1=L_{\mathcal{K}}^2\left( \mathbb{R}%
_{+}\right) $ the functional Hilbert space corresponding to the minimal
dilation (\ref{4.0}) sub-space $\mathcal{K}=\mathcal{K}_1$ for the CP map $%
\boldsymbol{\varphi }$, given by the orthoprojector $J=J_1$ in the space $%
\mathcal{K}_{\circ }$ of the canonical dilation, and $\mathfrak{K}_0$ its
orthogonal compliment, corresponding to $\mathcal{K}_0=J_0\mathcal{K}_{\circ
}$, where $J_0=I-J_1$. Representing $\Gamma \left( \mathfrak{K}_0\oplus 
\mathfrak{K}_1\right) $ as $\Gamma \left( \mathfrak{K}_0\right) \otimes
\Gamma \left( \mathfrak{K}_1\right) $, let us denote by $I_t$ the survival
orthoprojectors 
\begin{equation*}
I_t\chi \left( \sigma ^0,\sigma ^1\right) =\delta _\emptyset \left( \sigma
_t^0\right) \chi \left( \sigma ^0,\sigma ^1\right) ,\qquad \sigma _t=\sigma
\cap [0,t),
\end{equation*}
where $\chi \left( \sigma ^0,\sigma ^1\right) =\chi \left( \sigma ^0\sqcup
\sigma ^1\right) \in \mathcal{K}^{\otimes \left| \sigma ^0\right| }\otimes 
\mathcal{K}^{\otimes \left| \sigma ^1\right| }$ is the set function,
representing a $\chi \in \Gamma \left( \mathfrak{K}_0\oplus \mathfrak{K}%
_1\right) $. The decreasing family $\left( I_t\right) _{t>0}$ defines the
decay orthoprojectors $E_t=I-I_t$ in $\Gamma \left( \mathfrak{K}_{\circ
}\right) $ satisfying the quantum stochastic equation $\mathrm{d}%
E_t=E_tJ_0\cdot \mathrm{d}\Lambda _{\circ }^{\circ }$ with $E_0=0$, and $%
\Lambda _{\circ }^{\circ }$ is the number integrator in the Fock space $%
\Gamma \left( \mathfrak{K}_{\circ }\right) $ over $\mathfrak{K}_{\circ }=%
\mathfrak{K}_0\oplus \mathfrak{K}_1$. Then one easily find that the minimal
CP flow (\ref{3.8'}) can be represented as $\phi _t\left( B\right)
=V_t^{*}\left( I_t\otimes B\right) V_t $.

We may also construct by iteration the minimal quantum stochastic $*$-flow $%
\tau =\left( \tau _t\right) $ over the operator algebra $\mathcal{A}$,
resolving the quantum Langevin equation by its iteration. It is unique if
normaliesed $\tau _t\left( I\right) =I_t$ to the solution of the Langevin
equation for $A=I$, and it is $*$-multiplicative as in \cite{18} due to the
differential $*$-multiplicativity of $\mathbf{j}$. Then the properties (\ref%
{4.9}) follow from the definition of the operators $\widehat{V}_t$, and can
be checked recurrently by use of (\ref{4.6}) and (\ref{4.7}).

The cocycle $U=\left( U_t\right) $ is constructed to give the unitary
solution to the HP equation (\ref{1.10}), which always exists due to
differential unitarity relations. If the solution is unique as in the case
of all bounded coefficients $K$ and $L$, it can be represented in the form
of the stochastic multiple integral of the chronologically ordered products
of the coefficients of the quantum differential equation under the
integrability conditions given in \cite{Be97}.

If $K+K^{\dagger }\geq \varphi \left( I\right) $, the HP unitarity condition 
\cite{16} is satisfied for the canonical choice $L_{+}^{\circ }=L^{\circ }$
and arbitrary isometric operator $J_{\bullet }^{\circ }$, $J_{\circ
}^{\bullet }J_{\bullet }^{\circ }=I_{\bullet }^{\bullet }$ with $K_{\bullet
}^{-}=L_{\circ }J_{\bullet }^{\circ }$, $K_{\circ }^{-}=L_{\circ }J_{\circ
}^{\circ }$, $J_{\circ }^{\circ }=I_{\circ }^{\circ }-J_{\bullet }^{\circ
}J_{\circ }^{\bullet }$; if $\boldsymbol{K}+\boldsymbol{K}^{\dagger }\geq 
\boldsymbol{\varphi }\left( I\right) $, in addition we make the choice $%
J_{\bullet }^{\circ }=L_{\bullet }^{\circ }$ from the canonical dilation,
and so $K_{\bullet }^{-}=L_{\circ }L_{\bullet }^{\circ }=K_{\bullet }$,
where $L_{\circ }^{*}=L^{\circ }$, $J_{\circ }^{\bullet }=J_{\bullet
}^{\circ *}$. In the first, subfiltering case$\ \lambda \left( I\right) \leq
0$ such a choice gives the coincidence $U_t\left( \delta _\emptyset \otimes
\psi _0\right) =V_t\psi _0$ of the stochastic multiple integrals for any
initial vacuum $\psi _0=\eta \otimes \delta _\emptyset $, $\eta \in \mathcal{%
D}$, and therefore $\left\| V_t\psi _0\right\| =\left\| \psi _0\right\| $.
Thus $M_t=V_t^{*}V_t$ is a martingale, and the equation (\ref{4.11}) is
satisfied for any initial $\psi _0$. In the second, contractive case $%
\boldsymbol{\lambda }\left( I\right) \leq 0$ the canonical choice gives $%
U_t\left( \delta _\emptyset \otimes \psi \right) =V_t\psi $ and therefore $%
\left\| V_t\psi \right\| =\left\| \psi \right\| $ for any $\psi \in \mathcal{%
D}\otimes \mathfrak{F}$. Thus $V_t^{*}V_t=I$, and the equation (\ref{4.11})
is satisfied for any state $\psi $. 
\endproof%

\end{document}